\documentclass[a4paper,reqno,11pt]{amsart}

\usepackage{mathtools}
\usepackage{amsfonts}
\usepackage{amssymb}
\usepackage{graphicx,subcaption}
\usepackage[shortlabels]{enumitem}
\usepackage{csquotes}
\usepackage{microtype}
\usepackage{xcolor}
\usepackage{hyperref}
\usepackage{url}
\usepackage{cases}

\usepackage[margin=1in]{geometry}

\newtheorem{theorem}{Theorem}[section]
\newtheorem{lemma}[theorem]{Lemma}
\newtheorem{proposition}[theorem]{Proposition}
\newtheorem{corollary}[theorem]{Corollary}

\theoremstyle{definition}
\newtheorem{definition}[theorem]{Definition}

\theoremstyle{remark}

\newtheorem{remark}[theorem]{Remark}

\numberwithin{equation}{section}
\numberwithin{table}{section}
\numberwithin{figure}{section}

\newcommand\DN{\mathrm{DN}}
\renewcommand\i{\mathrm{i}}

\newcommand{\HGc}{H^{1/2}_\chi(\Gamma)}
\newcommand{\HmGc}{H^{-1/2}_\chi(\Gamma)}
\newcommand{\p}{\partial}
\newcommand{\gD}{\gamma_{_D}^\chi}
\newcommand{\gN}{\gamma_{_N}^\chi}
\newcommand{\cC}{\mathcal C}
\newcommand{\cD}{\mathcal D}
\newcommand{\cE}{\mathcal E}
\newcommand{\bbR}{\mathbb R}

\DeclareMathOperator{\Mor}{Mor}

\DeclareMathOperator{\Int}{Int}

\DeclareMathOperator{\dom}{dom}
\DeclareMathOperator{\spn}{span}
\let\phi=\varphi
\let\epsilon=\varepsilon

\title[Computing nodal deficiency with a refined Dirichlet-to-Neumann map]{Computing
  nodal deficiency with a refined Dirichlet-to-Neumann map}

\author{G. Berkolaiko}
\address{Department of Mathematics, Texas A\&M University, College Station, TX 77843-3368, USA}
\email{gberkolaiko@tamu.edu}

\author{G. Cox}
\address{Department of Mathematics and Statistics, Memorial University of Newfoundland, St. John's, NL A1C 5S7, Canada}
\email{gcox@mun.ca}

\author{B. Helffer}
\address{Laboratoire de Math\'{e}matiques Jean Leray,
Universit\'{e} de Nantes, 44 322 Nantes, France}
\email{Bernard.Helffer@univ-nantes.fr}

\author{M.P. Sundqvist}
\address{Lund University, Department of Mathematical Sciences, Box 118, 221 00 Lund, Sweden}
\email{mikael.persson\_sundqvist@math.lth.se }


\begin{document}

\begin{abstract}
Recent work of the authors and their collaborators has uncovered fundamental connections between the Dirichlet-to-Neumann map, the spectral flow of a certain family of self-adjoint operators, and the nodal deficiency of a Laplacian eigenfunction (or an analogous deficiency associated to a non-bipartite equipartition). Using a refined construction of the Dirichlet-to-Neumann map, we strengthen all of these results, in particular getting improved bounds on the nodal deficiency of degenerate eigenfunctions. Our framework is very general, allowing for non-bipartite partitions, non-simple eigenvalues, and non-smooth nodal sets. Consequently, the results can be used in the general study of spectral minimal partitions, not just nodal partitions of generic Laplacian eigenfunctions.
\end{abstract}

\maketitle


\section{Introduction}
Let $\Omega \subset \bbR^2$ be an open, bounded set, with piecewise $C^2$ boundary, and suppose $\phi_*$ is an eigenfunction of the Dirichlet Laplacian $-\Delta$ on $\Omega$, with eigenvalue $\lambda_*$. We denote by $\Gamma$ the nodal set of $\phi_*$,
\[
	\Gamma = \overline{\{x\in\Omega : \phi_*(x)=0\}},
\]
and by $k(\phi_*)$  the number of nodal domains of $\phi_*$, i.e.\ the number of connected components of the set $\{x\in\Omega : \phi_*(x)\neq 0\}$. We also let $\ell(\phi_*) = \min\{ m : \lambda_m = \lambda_*\}$ denote the minimal label of the eigenvalue $\lambda_*$, where $\lambda_1 < \lambda_2 \leq \lambda_3 \leq \cdots$ are the ordered Dirichlet eigenvalues of $\Omega$, repeated according to their multiplicity. The Courant nodal domain theorem states that $k(\phi_*) \leq \ell(\phi_*)$, or equivalently, that the nodal deficiency $\delta(\phi_*) := \ell(\phi_*) - k(\phi_*)$ is nonnegative.

Letting $\DN(\Gamma,\lambda_*)$ denote the two-sided Dirichlet-to-Neumann map on $\Gamma$, which will be defined below, we now state a special case of our main result.

\begin{theorem}
\label{thm:nodal}
The eigenfunction $\phi_*$ has nodal deficiency
\begin{equation}
\label{eq:nodal}
	\delta(\phi_*) = \Mor  \DN(\Gamma,\lambda_*),
\end{equation}
and the corresponding eigenvalue $\lambda_*$ has multiplicity
\begin{equation}
\label{eq:mult1}
	\dim\ker(\Delta + \lambda_*) = \dim\ker \DN(\Gamma,\lambda_*) + 1.
\end{equation}
\end{theorem}

The symbol $\Mor$ denotes the \emph{Morse index}, i.e.\ the number of negative eigenvalues of the operator $\DN(\Gamma,\lambda_*)$, which is self-adjoint and lower semi-bounded.
A similar formula for the nodal deficiency appeared in~\cite{CJM}; see also~\cite{BeCoMa}. The version of the Dirichlet-to-Neumann map appearing in the above theorem is more involved than the one used in~\cite{BeCoMa,CJM}, but consequently gives us a stronger result, as we now explain.

We denote the nodal domains of $\phi_*$ by $D_1, \ldots, D_k$.  When defining the Dirichlet-to-Neumann map, one must take into account that $\lambda_*$ is a Dirichlet eigenvalue on each $D_i$. Introducing the notation $\Gamma_i = \overline{\p D_i \cap \Omega}$, we define the closed subspace
\begin{equation}
\label{Sdef1}
	S = \left\{g\in L^2(\Gamma) : \int_{\Gamma_i} g_i \, \frac{\p \phi_{*,i}}{\p \nu_i}  = 0, \ i=1,\dots, k\right\}
\end{equation}
of $L^2(\Gamma)$, where $g_i$ denotes the restriction of $g$ to $\Gamma_i$, $\phi_{*,i}$ is the restriction of $\phi_{*}$ to $D_i$, and $\nu_i$ is the outward unit normal to $D_i$. For sufficiently smooth $g \in S$, each boundary value problem
\begin{equation}\label{eq:Di0}
	\begin{cases}
	-\Delta u_i = \lambda_* u_i &\text{in $D_i$},\\
		u_i = g_i & \text{on $\partial D_i\cap \Omega$},\\
		u_i=0 & \text{on $\partial D_i\cap\partial\Omega$},
	\end{cases}
\end{equation}
has a solution $u_i^g$.  Defining a function $\gamma_{_N} u^g$ on $\Gamma$ by
\begin{equation}
\label{2nudef0}
	\gamma_{_N} u^g\big|_{\Gamma_i\cap\Gamma_j} = \frac{\p u_i^g}{\p\nu_i} + \frac{\p u_j^g}{\p\nu_j}
\end{equation}
for all $i \neq j$, we let
\begin{equation}
\label{DN1}
	\DN(\Gamma,\lambda_*)g = \Pi_S \big(\gamma_{_N} u^g\big),
\end{equation}
where $\Pi_S$ denotes the $L^2(\Gamma)$-orthogonal projection onto the
subspace $S$.

The solution to the problem \eqref{eq:Di0} is non-unique, but the
choice of particular solution $u_i^g$ is irrelevant for the definition
on account of the projection in~\eqref{DN1}. In
Theorem~\ref{thm:Aform} we use this freedom to give an equivalent
formulation of the Dirichlet-to-Neumann map that does not involve
$\Pi_S$.

The earlier works~\cite{BeCoMa,CJM} avoided the difficulty of defining
the Dirichlet-to-Neumann map at a Dirichlet eigenvalue by evaluating
the quantities in Theorem~\ref{thm:nodal} at $\lambda_* + \epsilon$,
with a small positive $\epsilon$.  The resulting expression for the
nodal deficiency was
\begin{equation}
\label{eq:nodal2}
	\delta(\phi_*) = \Mor  \DN(\Gamma,\lambda_*+\epsilon) + 1 - \dim\ker(\Delta + \lambda_*).
\end{equation}
Unlike~\eqref{eq:nodal}, which immediately implies $\delta(\phi_*) \geq 0$, the equality~\eqref{eq:nodal2} only yields the same conclusion if we know that $\lambda_*$ is simple, or have additional information about the spectrum of $\DN(\Gamma,\lambda_*+\epsilon)$. Therefore, we obtain a more useful result by computing the Dirichlet-to-Neumann map at $\lambda_*$ instead of $\lambda_*+\epsilon$.

An even stronger motivation for eliminating the $\epsilon$-perturbation
is that the unperturbed operator $\DN(\Gamma,\lambda_*)$ appears
naturally as the Hessian of the energy functional 
on the space of generic equipartitions~\cite{BCCM}. The minima of this functional are  
spectral minimal partitions, as defined in~\cite{HHOT}, which are
often non-bipartite (unlike the decompositions of $\Omega$ into nodal
domains of an eigenfunction $\phi_*$, mentioned above).
One of the simplest examples of a non-bipartite partition is the
so-called Mercedes star partition, which is an (unproven but natural)
candidate for the minimal $3$-partition of the disk;
see~\cite{bonnaillie2015nodal} and references therein.  The main
result of this paper, Theorem~\ref{thm:general}, is a generalization
of Theorem~\ref{thm:nodal} to partitions that are not necessarily
bipartite, but have certain criticality properties that make
them prime candidates for being minimal.

We first recall\footnote{Here we are following the convention of \cite{HeSu};
  in \cite{bonnaillie2015nodal,HHOT} such a $\cD$ is called a
  \emph{strong partition}.} that a \emph{$k$-partition of $\Omega$} is
a family $\cD = \{D_i\}_{i=1}^k$ of mutually disjoint, open, connected
subsets of $\Omega$, with
$\overline\Omega = \overline{D_1 \cup \cdots \cup D_k}$.
We say that the subdomains $D_i$ and $D_j$  are \emph{neighbors}
if $\Int (\overline {D_i\cup
  D_j}) \neq D_i \cup D_j$.
We also recall that $\mathcal D$ is \emph{bipartite} if we can color
the partition with two colors in such a way that any two neighbors
have different colors.
Defining the
\emph{boundary set} of the partition to be
\begin{equation}
	\Gamma := \overline{\bigcup_i (\p D_i \cap \Omega)},
\end{equation}
we next impose a suitable regularity assumption on $\cD$.

\begin{definition}
\label{def:regular}
A partition $\cD$ is said to be \emph{weakly regular} if its boundary set $\Gamma$ satisfies:
\begin{enumerate}[(i)]
\item 
Except for finitely many critical points $ \{x_\ell\}\subset \Gamma\cap\Omega$,
 $\Gamma $ is locally diffeomorphic to a regular curve. In a
neighborhood of each  $ x_\ell$, $\Gamma$ is a union
of $\nu_\ell\geq 3$ smooth half-curves with one end at $ x_\ell$.
\item $ \Gamma\cap\partial\Omega$ consists of a 
finite set of boundary points $ \{z_m\}$. In a neighborhood of each 
 $ z_m$,  $ \Gamma$ is a union of $ \rho_m$
distinct smooth half-curves with one end at~$ z_m$.
\item The half-curves meeting at each $x_\ell$ and $z_m$ are pairwise transversal to one another, and to $\p\Omega$.
\end{enumerate}
\end{definition}

\begin{figure}[htb]
  \centering
  \includegraphics[page=7]{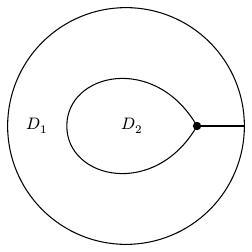}
  \caption{A partition of a disk that is not two-sided, since $D_1$ neighbors itself.}
  \label{fig:selfneighbor}
\end{figure}

The subdomains are only allowed to have corners at points where at
least three subdomains meet, or on $\partial\Omega$. However, the
definition still allows for partitions where a subdomain $D_i$
is a neighbor of itself, as shown in
Figure~\ref{fig:selfneighbor}. To rule out such examples, we say
that a partition $\cD$ is \emph{two-sided}\footnote{In \cite{bonnaillie2015nodal} such partitions are said to be \emph{nice}. We prefer the term two-sided, as it conveys the fact that each smooth 
component of $\Gamma$ is contained in the boundary of two distinct subdomains.} if $\Int(\overline D_i) = D_i$
for each $i$. 
For the rest of the paper we will only consider two-sided, weakly regular
partitions. This is a reasonable hypothesis, as it is satisfied by 
nodal partitions, and more generally by spectral minimal partitions \cite{Bers1955,HHOT}. 

For a two-sided, weakly regular partition each $D_i$ is a Lipschitz domain, so 
we can define trace operators and solve boundary value problems in a standard way. Without the 
transversality condition (iii) the $D_i$ may have cusps, 
and the analysis becomes much more difficult; see, for instance \cite{BCE}. If the partition is not two-sided, then some $D_i$ lies on both sides of its boundary. In this 
case it is possible to define separate trace operators on each side of the common 
boundary; we do not to this here, but refer to \cite[Section~1.7]{Gr} for an example 
of this construction.

To extend the notion of a ``nodal partition'' to partitions that are
not necessarily bipartite, it is convenient to introduce a generalization 
of the Laplacian. The construction involves a choice of signed weight
functions, which will also be used to define a generalized two-sided
Dirichlet-to-Neumann map on the partition boundary set.  

\begin{definition}
  \label{def:weight}
  Given a two-sided, weakly regular partition $\cD=\{D_i\}$, let
  \begin{equation}
    \label{eq:Gamma-i-def}
    \Gamma_i := \overline{\p D_i \cap \Omega}.  
  \end{equation}
  We say that functions $\chi_i \colon \Gamma_i \to \{\pm1\}$ are
  \emph{valid weights} if they are constructed as follows. Given an 
  orientation of each $\p D_i$, and an orientation of
  each smooth component of $\Gamma$, we define $\chi_i$ on each 
  smooth component of $\Gamma_i$ to be $+1$ if the
  orientation of $\p D_i$ agrees with the orientation of the 
  corresponding smooth component of $\Gamma$, and
  equal to $-1$ otherwise.
\end{definition}

\begin{figure}[htp]
  \centering
  \subcaptionbox{\label{figa}}{\includegraphics[page=1,scale=0.8]{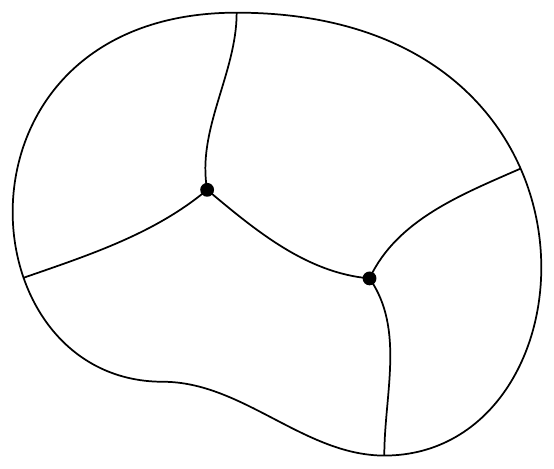}}
  \hskip 1cm
  \subcaptionbox{\label{figb}}{\includegraphics[page=2,scale=0.8]{figures3.pdf}}
  \hskip 1cm
  \subcaptionbox{\label{figc}}{\includegraphics[page=3,scale=0.8]{figures3.pdf}}
  \hskip 1cm
  \subcaptionbox{\label{figd}}{\includegraphics[page=4,scale=0.8]{figures3.pdf}}
  \hskip 1cm
  \subcaptionbox{\label{fige}}{\includegraphics[page=5,scale=0.8]{figures3.pdf}}
  \hskip 1cm
  \subcaptionbox{\label{figf}}{\includegraphics[page=6,scale=0.8]{figures3.pdf}}
  \hskip 1cm
  \caption{ A partition \subref{figa}, with a choice of orientation
    for the boundary $\partial D_i$ of each subdomain \subref{figb},
    and an orientation of each smooth part of the boundary set
    $\Gamma$ \subref{figc}.  In \subref{figd} we show the
    corresponding valid weights $\chi_i$, and in \subref{fige} we show
    the resulting cut, as described in Appendix~\ref{sec:weights}. In
    \subref{figf} we display a non-valid choice of weights, i.e.\
    functions $\chi_i \colon \Gamma_i \to \{\pm1\}$ that are not induced by any
    choice of orientations.  }
  \label{fig:partitionwithcutsandsigns}
\end{figure}

In Figure~\ref{fig:partitionwithcutsandsigns} we illustrate the construction of a valid set of weights, and also give an example of a non-valid choice of weights. Note that $\chi_i$ is constant on each smooth segment of $\Gamma_i$; 
the value at the corner
points is irrelevant. According to Definition~\ref{def:weight}, there
are two ways $\chi_i$ can change sign on $\Gamma_i$: 1)
it can change sign at a corner; or 2) it can take different
signs on different connected components. It is easily shown that a partition is
bipartite if and only if the weights $\chi_i\equiv 1$ are valid, cf. \cite[Lemma~9]{BCCM}, 
and so non-constant weights are essential for the study of non-bipartite partitions.

\begin{remark}
\label{rem:coorientation}
An equivalent construction of valid weights can be given in terms of a co-orientation of each $\p D_i$ 
and each smooth component of $\Gamma$. Along each $\p D_i$ we choose a vector field 
$V_i$ that is equal to either $\nu_i$ or $-\nu_i$. Choosing a vector
field $V$ that is a smooth unit normal to each smooth component of
$\Gamma$, we set $\chi_i = V \cdot V_i$. A special case of this
construction appeared in \cite{BCCM}, where $V_i$ was chosen to
be the outward unit normal $\nu_i$, in which case $\chi_i = - \chi_j$
whenever $D_i$ and $D_j$ are neighbors. The extra flexibility in the present 
construction will be useful below, in our discussion of $\chi$-nodality.
\end{remark}

Valid weights have a natural geometric interpretation in terms of the cutting construction in~\cite[Section~4]{HeSu}, where one removes a portion $\Gamma^*$ of the nodal set from the domain $\Omega$ in such a way that the resulting partition of $\Omega \setminus \Gamma^*$ is bipartite; see Appendix~\ref{sec:weights} for details.

We now introduce a weighted version of the Laplacian, $-\Delta^\chi$,
corresponding to the bilinear form defined on the domain
$\dom(t^\chi)$ consisting of $u \in L^2(\Omega)$ such that
\begin{numcases}{}
	u_i := u\big|_{D_i} \in H^1(D_i), \label{eq:domt_H1} \\
	u_i = 0 \text{ on } \p D_i \cap \p\Omega, \label{eq:domt_Dext} \\
	\chi_i u_i = \chi_j u_j \text{ on } \Gamma_i \cap \Gamma_j \text{
  for all } i,j =1,\ldots, k, \label{eq:domt_Dint}
\end{numcases}
and given by
\begin{equation}
  \label{formt}
  t^\chi(u,v) = \sum_{i=1}^k \int_{D_i} \nabla u_i \cdot \nabla v_i.
\end{equation}
The Laplacians $\Delta^\chi$ for different valid weights will be shown
in Proposition~\ref{prop:Deltachi} to be unitarily equivalent.  As a
consequence, if the partition is bipartite, then $\Delta^\chi$ is unitarily
equivalent to the Dirichlet Laplacian on $\Omega$.  Furthermore, the
nodal sets of the eigenfunctions of $\Delta^\chi$ are independent of
$\chi$, justifying the following definition.

\begin{definition}
  \label{def:chinodal}
  A two-sided, weakly regular partition $\cD$ is said to be
  \emph{$\chi$-nodal} if it is the nodal partition for some
  eigenfunction of $\Delta^\chi$.  The \emph{defect} of a $\chi$-nodal
  $k$-partition is defined to be
  \begin{equation}
    \delta(\cD) = \ell(\cD) - k,
  \end{equation}
  where $\ell(\cD)$ denotes the minimal label of $\lambda_*$ in the
  spectrum of $-\Delta^\chi$.
\end{definition}

In Section~\ref{sec:chiSPCC} we will show that a partition is
$\chi$-nodal if and only if it satisfies the strong pair compatibility
condition \cite{HeSu}.

\begin{definition}
\label{def:SPCC1}
  A two-sided, weakly regular partition $\cD$ is said to satisfy the
  \emph{strong pair compatibility condition} (SPCC) if there exists a
  choice of positive ground states $\{u_i\}_{i=1}^k$
  for the Dirichlet Laplacians on $D_i$ such that, for any pair of neighbors
  $ D_i$ and $D_j$, the function $u_{ij}$ defined by
  \begin{equation}
    \label{eq:u_ij}
    u_{ij} \big|_{D_i} = u_i,
    \qquad
    u_{ij} \big|_{D_j} = -u_j,    
  \end{equation}
  is an eigenfunction of the Dirichlet Laplacian on
  $\Int (\overline {D_i\cup D_j})$.
\end{definition}

 We stress
that the choice of the ground states in the definition (which is
merely a choice of normalization on each $D_i$) is global\,---\,it 
can not change from one pair of neighbors to another.
This distinguishes SPCC from the \emph{weak pair compatibility 
condition}\footnote{In earlier papers, for instance \cite{HH:2005a}, WPCC is simply referred 
to as the \emph{pair compatibility condition} (PCC).} (WPCC) also appearing in
the literature; see Appendix~\ref{sec:weights}. 
It is immediate that nodal partitions satisfy the SPCC. We also mention that 
for a smooth partition, where the set 
$ \{x_\ell\}\subset \Gamma\cap\Omega$ 
of singular points is empty, the SPCC is equivalent to being 
a critical point of the energy functional on the set of
equipartitions; see \cite{BKS}.

Finally, we will define a $\chi$-weighted version of the two-sided
Dirichlet-to-Neumann map, denoted $\DN(\Gamma,\lambda_*,\chi)$.  The
full definition, given in Section~\ref{sec:def}, is rather delicate
because $\lambda_*$ is a Dirichlet eigenvalue and $\Gamma$ has
corners.  We just mention here that, similar to the Laplacian
$\Delta^\chi$, the Dirichlet-to-Neumann maps defined with different
valid $\{\chi_i\}$ are unitarily equivalent; the precise nature of the
equivalence is clarified in Theorem~\ref{thm:Aform}.  If each
$\chi_i$ is constant, $\DN(\Gamma,\lambda_*,\chi)$ reduces to the
operator $\DN(\Gamma,\lambda_*)$ already described in \eqref{DN1}.

The main result of this paper is the following. 

\begin{theorem}
  \label{thm:general}
A two-sided, weakly regular partition $\cD$ satisfies the SPCC if and only if it is $\chi$-nodal, in which case it has defect
  \begin{equation}
    \label{eq:def}
    \delta(\cD) = \Mor  \DN(\Gamma,\lambda_*,\chi),
  \end{equation}
  and the corresponding eigenvalue $\lambda_*$ of $-\Delta^\chi$ has multiplicity
  \begin{equation}
    \label{eq:mult2}
    \dim\ker(\Delta^\chi + \lambda_*) = \dim\ker \DN(\Gamma,\lambda_*,\chi) + 1.
  \end{equation}
  The quantities in \eqref{eq:def} and \eqref{eq:mult2} are
  independent of $\chi$.  In particular, different valid weights
  $\{\chi_i\}$ may be used in defining the Laplacian $-\Delta^\chi$
  and the Dirichlet-to-Neumann map $\DN(\Gamma,\lambda_*,\chi)$.
\end{theorem}

\begin{remark}
  Theorem~\ref{thm:general} contains Theorem~\ref{thm:nodal} as a
  special case, and hence is an improvement over the results
  of~\cite{BeCoMa,CJM}, as described above. Similarly, for
  non-bipartite partitions, it refines~\cite[Theorem~4.1]{HeSu}.
\end{remark}

\begin{remark}
 In higher dimensions the nodal sets of eigenfunctions can be more 
 complicated, and the analysis of corner domains is significantly more involved 
 (see, for instance \cite{Dauge}), so we restrict our attention to the planar case.
 The conclusion of
  Theorem~\ref{thm:general} immediately extends to higher dimensions
  if the nodal set $\cup \Gamma_i$ is a smoothly embedded
  hypersurface.
\end{remark}

\subsection*{Outline}
In Section~\ref{sec:chinodal} we give some preliminary analysis, describing Sobolev spaces on the boundary set $\Gamma$, weighted Dirichlet and Neumann traces, and the weighted Laplacian $\Delta^\chi$. We also show that a partition is $\chi$-nodal if and only if it satisfies the SPCC, and 
prove some delicate regularity results.
In Section~\ref{sec:def} we define the weighted, two-sided Dirichlet-to-Neumann operator $\DN(\Gamma,\lambda_*,\chi)$ and establish its fundamental properties. In Section~\ref{sec:flow} we prove Theorem~\ref{thm:general} by studying the spectral flow of an analytic family of self-adjoint operators. In Section~\ref{sec:circle} we illustrate our results by applying them to partitions of the circle.

In Appendix~\ref{sec:weights} we discuss the strong and weak pair compatibility conditions, and the connection between our $\chi$ weights and the cutting construction of \cite{HeSu}. Finally, in Appendix~\ref{sec:equivalent} we describe an alternate, more explicit construction of the canonical solution to a boundary value problem that arises in our construction of the Dirichlet-to-Neumann map.

\subsection*{Acknowledgments}
The authors thank Yaiza Canzani, Jeremy Marzuola and Peter Kuchment
for inspiring discussions about nodal partitions and
Dirichlet-to-Neumann maps, and the organizers of the \emph{Spectral
  Geometry in the Clouds} seminar (namely, Alexandre Girouard, Jean
Lagac\'e and Laura Monk), where the present collaboration was initiated.
G.B. acknowledges the support of NSF Grant DMS-1815075.
G.C. acknowledges the support of NSERC grant RGPIN-2017-04259.

\section{Preliminary analysis}
\label{sec:chinodal}

In this section we provide some background for our construction
of the Dirichlet-to-Neumann map, in particular defining Sobolev spaces
on the boundary set $\Gamma$, weighted Dirichlet and Neumann traces,
and the weighted Laplacian.  We also establish that SPCC is equivalent
to $\chi$-nodality.

\subsection{Sobolev spaces on the boundary set}
\label{sec:Sobolev}

Recall that $\Gamma_i = \overline{\p D_i \cap \Omega}$.
Since $|\chi_i| \equiv 1$ on
$\Gamma_i$, we have
\begin{align*}
  g \in L^2(\Gamma)
  &\quad \Longleftrightarrow \quad
  g_i := g\big|_{\Gamma_i} \in L^2(\Gamma_i) \ \text{ for each } i \\
  &\quad \Longleftrightarrow \quad
  \chi_i g_i \in L^2(\Gamma_i) \ \text{ for each } i.  
\end{align*}
The situation for $H^{1/2}$ is more complicated. If $\Gamma$
has intersections then it is not a Lipschitz manifold, and the space
$H^{1/2}(\Gamma)$ cannot be defined in the usual way; cf. \cite{McL}.  Moreover, on each
subdomain the conditions $g_i \in H^{1/2}(\Gamma_i)$ and
$\chi_i g_i \in H^{1/2}(\Gamma_i)$ need not be equivalent, due to the
possible discontinuities of $\chi_i$ at the corner points. We thus
define the space
\begin{equation}
\label{def:HGc}
	\HGc := \big\{ g \in L^2(\Gamma) : \cE_i (\chi_i g_i) \in H^{1/2}(\p D_i), \ i=1,\dots, k \big\},
\end{equation}
where $\cE_i \colon L^2(\Gamma_i) \to L^2(\p D_i)$ is the extension by zero to the rest of $\p D_i$, i.e.
\[
	\cE_i(\chi_i g_i) := \begin{cases} \chi_i g_i & \text{on } \Gamma_i, \\
	0  & \text{on } \p D_i \setminus \Gamma_i. \end{cases}
\]
The condition $\cE_i(\chi_i g_i) \in H^{1/2}(\p D_i)$ is more
restrictive than $\chi_i g_i \in H^{1/2}(\Gamma_i)$ if $\p D_i \cap
\p\Omega \neq \varnothing$. For instance, if $\chi_i g_i$ is a nonzero
constant on $\Gamma_i$, its extension by zero will not be an element
of $H^{1/2}(\p D_i)$. A necessary and sufficient condition for
$\cE_i(\chi_i g_i) \in H^{1/2}(\p D_i)$ will be recalled below, in
Lemma~\ref{lem:glue}. We define the norm
\begin{equation}
  \label{norm:HGc}
  \|g\|^2_{\HGc} := \sum_{i=1}^k \left\| \cE_i(\chi_i g_i)\right\|^2_{H^{1/2}(\p D_i)},
\end{equation}
and let $\HmGc$ denote the dual space to $\HGc$.

We next define a weighted Dirichlet trace (i.e.\ restriction to the
nodal set) operator. A natural domain for this operator is the set
$\dom(t^\chi)$ that was defined above
in~\eqref{eq:domt_H1}--\eqref{eq:domt_Dint}, equipped with the norm
$\|u\|^2_{\dom(t^\chi)} = \sum_i \| u_i \|^2_{H^1(D_i)}$.

\begin{lemma}
\label{lem:Dtrace}
The trace map
\begin{equation}
	\gD \colon \dom(t^\chi) \longrightarrow \HGc
\end{equation}
defined by $(\gD u)\big|_{\Gamma_i} = \chi_i u_i\big|_{\Gamma_i}$ is bounded, and has a bounded right inverse.
\end{lemma}

\begin{proof}
For each $D_i$ there is a bounded trace operator $H^1(D_i) \to H^{1/2}(\p D_i)$. We thus let $(\gD u)\big|_{\Gamma_i} = \chi_i u_i\big|_{\Gamma_i}$ for each $i$; the condition $\chi_i u_i = \chi_j u_j$ guarantees that $\gD u$ is a well-defined function for any $u \in \dom(t^\chi)$. Moreover, for each $i$ we have $\chi_i (\gD u)|_{\Gamma_i} = u_i\big|_{\Gamma_i}$, and hence
\[
	\cE_i \bigl(\chi_i (\gD u)\big|_{\Gamma_i}\bigr) = u_i\big|_{\p D_i}
\]
because $u_i = 0$ on $\p D_i \cap \p\Omega$. Since $u_i\big|_{\p D_i} \in H^{1/2}(\p D_i)$, it follows from~\eqref{def:HGc} that $\gD u \in \HGc$, with
\[
	\bigl\|\gD u\bigr\|^2_{\HGc} = \sum_{i=1}^k \bigl\| u_i|_{{\p D_i}} \bigr\|^2_{H^{1/2}(\p D_i)}
	\leq C \sum_{i=1}^k \| u_i \|^2_{H^1(D_i)},
\]
as was to be shown.

To construct a right inverse, we first recall that for each $i$ the trace map $H^1(D_i) \to H^{1/2}(\p D_i)$ has a bounded right inverse, $\Upsilon_i \colon H^{1/2}(\p D_i) \to H^1(D_i)$. Let $g \in \HGc$, so that $\cE_i(\chi_i g_i) \in H^{1/2}(\p D_i)$, and define $u_i = \Upsilon_i\big(\cE_i(\chi_i g_i)\big) \in H^1(D_i)$. The corresponding function $u \in L^2(\Omega)$, defined by $u\big|_{D_i} = u_i$ for each $i$, is contained in $\dom(t^\chi)$, since
\[
	\chi_i u_i\big|_{\Gamma_i \cap \Gamma_j} = \chi_i(\chi_i g_i) = g_i = g_j = \chi_j u_j\big|_{\Gamma_i \cap \Gamma_j}
\]
for all $i,j$. Moreover, we have
\[
	\|u\|^2_{\dom(t^\chi)} = \sum_i \|u_i\|^2_{H^1(D_i)} \leq C \sum_i \big\|\cE_i(\chi_i g_i)\big\|^2_{H^{1/2}(\p D_i)} = C \|g\|^2_{\HGc},
\]
and so $\Upsilon g = u$ defines a bounded right inverse $\Upsilon \colon \HGc \to \dom(t^\chi)$.
\end{proof}

We next define a weighted, two-sided version of the normal derivative that will appear naturally in our construction of the Dirichlet-to-Neumann map.

\begin{lemma}
\label{lemma:Ntrace}
If $u \in L^2(\Omega)$, with $u_i \in H^1(D_i)$ and $\Delta u_i \in L^2(D_i)$ for each $i$, then there exists a unique $\gN u \in \HmGc$ such that
\begin{equation}
\label{Ntrace}
	\left<\gN u, \gD v \right>_* = \sum_{i=1}^k \int_{D_i} \big(\nabla u_i \cdot \nabla v_i + (\Delta u_i)v_i \big)
\end{equation}
for all $v \in \dom(t^\chi)$, where $\left<\cdot,\cdot\right>_*$ denotes the dual pairing between $\HGc$ and $\HmGc$.
\end{lemma}

Note that the definition of $\gN u$ does not require any consistency conditions on the boundary values of $u$ along $\Gamma$. That is, we do not require $u \in \dom(t^\chi)$.

\begin{proof}
The construction is almost identical to that of~\cite[Lemma~4.3]{McL}. Letting $\Upsilon \colon \HGc \to \dom(t^\chi)$ denote a bounded right inverse to $\gD$, as in Lemma~\ref{lem:Dtrace}, we define $\gN u \in \HmGc$ by its action on arbitrary $g \in \HGc$, namely
\[
	\left<\gN u, g \right>_* := \sum_{i=1}^k \int_{D_i} \big(\nabla u_i \cdot \nabla (\Upsilon g)_i + (\Delta u_i)(\Upsilon g)_i \big).
\]
It is easily verified that this has all the required properties.
\end{proof}

\begin{remark}
\label{rem:Nreg}
If, in addition to the hypotheses of Lemma~\ref{lemma:Ntrace}, $u_i\big|_{\p D_i} \in H^1(\p D_i)$ for each $i$, then \cite[Theorem~4.24]{McL} implies 
$\p_{\nu_i} u_i \in L^2(\p D_i)$, and so Green's formula can be written as
\[
	\int_{D_i} \big(\nabla u_i \cdot \nabla v_i + (\Delta u_i)v_i \big) = \int_{\p D_i} \frac{\p u_i}{\p \nu_i} v_i = \int_{\Gamma_i} \left(\chi_i \frac{\p u_i}{\p \nu_i}\right) (\chi_i v_i).
\]
That is, the dual pairing of $\p_{\nu_i} u_i \in H^{-1/2}(\p D_i)$ and $v_i\big|_{\p D_i} \in H^{1/2}(\p D_i)$ is given by their $L^2(\p D_i)$ inner product.
Summing over $i$ and comparing with~\eqref{Ntrace}, we find that
\begin{equation}
  \label{gNfunction}
  \gN u \big|_{\Gamma_i \cap \Gamma_j}
  = \chi_i \frac{\p u_i}{\p \nu_i} + \chi_j \frac{\p u_j}{\p \nu_j} \in L^2(\Gamma_i \cap \Gamma_j)
\end{equation}
for any neighbors $D_i$ and $D_j$.
\end{remark}

\subsection{The sign-weighted Laplacian}

In this section we describe the self-adjoint operator $\Delta^\chi$
and its dependence on $\chi$.

\begin{definition}
  \label{def:weight_equivalences}
  We say that two sets of valid weights $\{\chi_i\}$ and
  $\{\tilde\chi_i\}$ are \emph{edge equivalent} if for each $i \neq j$
  we have $\tilde\chi_i \tilde\chi_j = \chi_i \chi_j$ on
  $\Gamma_i \cap \Gamma_j$, and \emph{domain equivalent} if for each
  $i$ we have either $\tilde\chi_i \equiv \chi_i$ or
  $\tilde\chi_i \equiv -\chi_i$.
\end{definition}

\begin{figure}[htp]
   \centering
   \subcaptionbox{\label{omegasigns1}}{\includegraphics[page=2,scale=0.8]{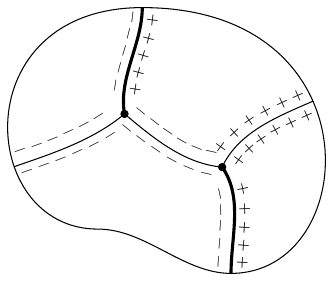}}
   \hskip 1cm
   \subcaptionbox{\label{omegasigns2}}{\includegraphics[page=3,scale=0.8]{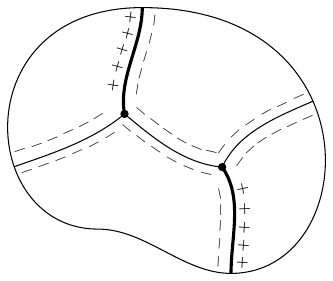}}
   \caption{The weights in \subref{omegasigns1} and \subref{omegasigns2} are edge equivalent but not domain equivalent.}
   \label{fig:edgeequivalent}
 \end{figure}

In Figure~\ref{fig:edgeequivalent}, the weights in~\subref{omegasigns1} and~\subref{omegasigns2} are edge equivalent.

\begin{remark}
  \label{rem:chi_all_equiv}
  In terms of Definition~\ref{def:weight}, edge
  equivalence corresponds to only changing the orientations of the
  smooth components of $\Gamma$, while domain equivalence corresponds
  to only changing the orientations of the $\p D_i$.  It is thus clear that
  for any valid sets of weights $\chi$ and $\tilde{\chi}$ there is a
  valid weight $\hat\chi$ such that $\chi$ is edge equivalent to
  $\hat\chi$ and $\hat\chi$ is domain equivalent to $\tilde\chi$.
\end{remark}

Recall that $-\Delta^\chi$ corresponds to the bilinear form $t^\chi$
defined in~\eqref{formt}, with $\dom(t^\chi)$ given
by~\eqref{eq:domt_H1}--\eqref{eq:domt_Dint}. The following proposition
summarizes its basic properties.

\begin{proposition}
  \label{prop:Deltachi}
  If $\cD$ is a two-sided, weakly regular partition and $\{\chi_i\}$ are
  valid weights, then $\Delta^\chi$ is a self-adjoint operator on
  $L^2(\Omega)$, with domain
  \begin{equation}
    \label{domD}
    \dom(\Delta^\chi) = \big\{ u \in \dom(t^\chi) \colon
    \Delta u_i \in L^2(D_i) \text{ for each } i \text{ and } \gN u = 0  \big\}.
  \end{equation}
  For any other set of valid weights $\{\tilde\chi_i\}$ we have:
  \begin{enumerate}
  \item If $\chi$ and $\tilde\chi$ are edge equivalent, then
    $\Delta^\chi = \Delta^{\tilde\chi} $;
  \item If $\chi$ and $\tilde\chi$ are domain equivalent, then
    $\Delta^\chi$ is unitarily equivalent to $\Delta^{\tilde\chi}$.
  \end{enumerate}
  Consequently, $\Delta^\chi$ and $\Delta^{\tilde\chi}$ are unitarily equivalent 
  for any choices of valid weights, and so the property of being $\chi$-nodal is independent of
  the choice of a valid $\chi$.
\end{proposition}

\begin{proof}
It is easily seen that $t^\chi$ is a closed, semi-bounded bilinear form, with dense domain in $L^2(\Omega)$. It thus generates a semi-bounded self-adjoint operator, which we denote $-\Delta^\chi$, with domain
\begin{align}
\label{domD2}
\begin{split}
	\dom(\Delta^\chi) = \big\{ u \in \dom(t^\chi)& : \text{there exists } f\in L^2(\Omega) \text{ such that } \\ &t^\chi(u,v) = \left<f,v\right>_{L^2(\Omega)}
	\text{ for all } v \in \dom(t^\chi)  \big\}.
\end{split}
\end{align}
For any such $u$ we have $-\Delta^\chi u = f$.

To prove~\eqref{domD}, we first assume that $u \in \dom(\Delta^\chi)$, as described in~\eqref{domD2}. If $v_i \in H^1_0(D_i)$ for some $i$, then its extension by zero to the rest of $\Omega$ is contained in $\dom(t^\chi)$. Denoting this extension by $v$, we get from~\eqref{formt} and~\eqref{domD2} that
\[
	\int_{D_i} fv_i = t^\chi(u,v) = \int_{D_i} \nabla u_i \cdot \nabla v_i.
\]
Since $v_i$ was an arbitrary function in $H^1_0(D_i)$, this means $\Delta u_i = f_i \in L^2(D_i)$ in a distributional sense. This holds for each $i$, so it follows from Lemma~\ref{lemma:Ntrace} that $\gN u$ is defined, and satisfies
\[
	\left<\gN u, \gD v \right>_* = t^\chi(u,v) -  \left<f,v\right>_{L^2(\Omega)} = 0
\]
for all $v \in \dom(t^\chi)$.  Since $\gD$ is surjective, this implies $\gN u = 0$.

On the other hand, suppose $u \in \dom(t^\chi)$ satisfies $\Delta u_i \in L^2(D_i)$ for each $i$ and $\gN u = 0$. Lemma~\ref{lemma:Ntrace} then implies
\[
	t^\chi(u,v) = - \sum_{i=1}^k \int_{D_i} (\Delta u_i)v_i = \left<f,v\right>_{L^2(\Omega)}
\]
for all $v \in \dom(t^\chi)$, where $f \in L^2(\Omega)$ is defined by $f_i = -\Delta u_i$ for each $i$. Using~\eqref{domD2}, this gives $u \in \dom(\Delta^\chi)$ and completes the proof.

Finally, we describe the dependence of the operator $\Delta^\chi$ on
the weights $\{\chi_i\}$.  The first claim follows immediately from
the definitions. If $\tilde\chi_i \tilde\chi_j = \chi_i \chi_j$ on
$\Gamma_i \cap \Gamma_j$, then $\chi_i u_i = \chi_j u_j$ is equivalent
to $\tilde\chi_i u_i = \tilde\chi_j u_j$, hence
$\dom(t^{\tilde\chi}) = \dom(t^\chi)$ and the result follows.

For the second claim, consider the unitary map $U \colon L^2(\Omega) \to L^2(\Omega)$ defined by
\[
	(Uv)\big|_{D_i} = \begin{cases} v & \text{ if } \tilde\chi_i \equiv \chi_i, \\
	-v & \text{ if } \tilde\chi_i \equiv -\chi_i. \end{cases}
\]
This sends $\dom(t^\chi)$ to $\dom(t^{\tilde\chi})$, with $t^{\tilde\chi}(Uv,Uw) = t^\chi(v,w)$ for all $v,w \in \dom(t^\chi)$, which implies $\Delta^\chi = U^{-1} \Delta^{\tilde\chi} U$ and completes the proof.

These two equivalences combined with
Remark~\ref{rem:chi_all_equiv} shows that $\Delta^\chi$ is unitarily
equivalent to any other $\Delta^{\tilde\chi}$ with a valid
$\tilde\chi$.  The unitary map $U$ does not affect the nodal set,
therefore a partition $\cD$ is $\chi$-nodal either for all valid
choices of $\chi$ or for none.
\end{proof}

\subsection{Pair compatibility and $\chi$-nodal partitions}
\label{sec:chiSPCC}
Next, we discuss the connection between the strong pair
compatibility condition and the $\chi$-nodal condition.

\begin{proposition}
  \label{prop:chinodal_spcc}
  A two-sided, weakly regular partition $\cD$ is $\chi$-nodal if and only
  if it satisfies the SPCC.
\end{proposition}

\begin{proof}
  First suppose $\cD$ is $\chi$-nodal, so it is the nodal set of some
  eigenfunction $\phi_*$ of $\Delta^\chi$. Since $\phi_{*,i}\big|_{\p D_i} 
  = 0$ is contained in $H^1(\p D_i)$, we can use Remark~\ref{rem:Nreg} and 
  Proposition~\ref{prop:Deltachi} to get
  \begin{equation}
  \label{eq:transmission}
    0 = \gN \phi_* \big|_{\Gamma_i \cap \Gamma_j}
    = \chi_i \frac{\p \phi_{*,i}}{\p \nu_i} + \chi_j \frac{\p \phi_{*,j}}{\p \nu_j}
  \end{equation}
for any neighbors $D_i$ and $D_j$.

  We now let $\eta_i = \operatorname{sgn} \phi_{*,i}$.  For each $D_i$
  the function $u_i := \eta_i \phi_{*,i}$ is a positive ground state, 
  and the
  transmission condition \eqref{eq:transmission} becomes
  \begin{equation}
    \label{eq:u-transmission}
    \eta_i\chi_i \frac{\p u_i}{\p \nu_i} + \eta_j\chi_j \frac{\p
      u_j}{\p \nu_j} = 0.
  \end{equation}
  Since $\p_{\nu_i} u_i$ and $\p_{\nu_j} u_j$ are both negative, we conclude that
  $\eta_i\chi_i = -\eta_j\chi_j$, yielding
  $\p_{\nu_i} u_i = \p_{\nu_j} u_j$ on $\Gamma_i \cap \Gamma_j$. It
  follows that $u_{ij}$, as defined in \eqref{eq:u_ij}, is a Dirichlet eigenfunction on $\Int (\overline {D_i\cup D_j})$,
  hence $\cD$ satisfies the SPCC.

  Conversely, suppose $\cD$ satisfies the SPCC. This means there exist
  positive ground states $u_i$ for the Dirichlet Laplacian on $D_i$
  such $u_{ij}$ is a Dirichlet eigenfunction on $\Int (\overline {D_i\cup D_j})$ whenever
  $D_i$ and $D_j$ are neighbors. This implies
  $\p_{\nu_i}u_i = \p_{\nu_j}u_j$ on $\Gamma_i \cap \Gamma_j$.
  
  Now define valid weights by choosing the same orientation for all $D_i$.
  This implies that $\chi_i = -\chi_j$ on $\Gamma_i \cap \Gamma_j$
  (the orientation of the segments of $\Gamma$ is irrelevant).  It
  follows that the function $u$ defined by $u\big|_{D_i} = u_i$
  satisfies the transmission condition \eqref{eq:transmission} and
  hence is an eigenfunction of $\Delta^\chi$.
\end{proof}

\begin{remark}
\label{rem:max}
  The $\chi$ weights chosen in the final step of the proof coincide with those used in \cite{BCCM}; 
  see Remark~\ref{rem:coorientation}. In Appendix~\ref{sec:weights} we will see that this $\chi$ corresponds to the so-called ``maximal cut'' of the boundary set.
\end{remark}

It is well known that nodal partitions and spectral minimal partitions have the \emph{equal angle property}:  at a singular point, the half-curves meet with equal angle \cite{HHOT}. This is also true of $\chi$-nodal partitions.

\begin{corollary}
\label{cor:equalangle}
If $\cD = \{D_i\}$ is $\chi$-nodal, then it satisfies the equal angle property.
\end{corollary}

\begin{proof}
Since $\cD$ satisfies the SPCC, the result follows from applying \cite[Theorem~2.6]{HHOT} 
to each pair of neighboring domains.
\end{proof}

It is an immediate consequence of the equal angle property that each $D_i$ has 
convex corners, a fact we will use in Proposition~\ref{prop:regDirichlet} to conclude 
$H^2$ regularity of Dirichlet eigenfunctions.

\subsection{Regularity properties of the Dirichlet kernel}
\label{sec:regularity}

Let $\Delta^\chi_\infty$
be the Laplacian in $\Omega$ with Dirichlet boundary conditions
imposed on $\partial\Omega\cup\Gamma$.  More precisely, it is the
Laplacian with the domain
\begin{align*}
  \dom(\Delta^\chi_\infty) &= \big\{ u \in \dom(t^\chi) :
  \Delta u_i \in L^2(D_i) \text{ for each } i \text{ and } \gD u = 0  \big\} \\
  &= \big\{ u \in L^2(\Omega) :
  u_i \in H^1_0(D_i) \text{ and } \Delta u_i \in L^2(D_i) \text{ for each } i  \big\}.
\end{align*}
The reason for the subscript $\infty$ will become apparent in
Section~\ref{sec:flow}.  For now, we would like to understand the
properties of the eigenspace of $-\Delta^\chi_\infty$ corresponding to
the eigenvalue $\lambda_*$.  The main result of this section is the
following.

\begin{proposition}
  \label{prop:regDirichlet}
  Let $\cD$ be a $\chi$-nodal $k$-partition and $\phi_*$ be the
  eigenfunction of $\Delta^\chi$ with boundary set $\Gamma$.  The subspace
  \begin{equation}
    \label{eq:Phi_space}
    \Phi := \ker\left(\Delta^\chi_\infty + \lambda_*\right)
  \end{equation}
  has the following properties:
  \begin{enumerate}
  \item $\dim \Phi = k$;
  \item $\ker \gN\big|_\Phi = \spn\{ \phi_* \}$;
  \item for any $\phi \in \Phi$,  $\phi\big|_{D_i} \in H^2(D_i) \cap H^1_0(D_i)$;
  \item $\gN(\Phi) \subset \HGc$.
  \end{enumerate}
\end{proposition}

\begin{proof}
  It follows immediately that for each $i$ the restriction
  $\phi_{*,i} \in H^1_0(D_i)$ of $\phi_*$ satisfies the eigenvalue
  equation $\Delta \phi_{*,i} + \lambda_* \phi_{i} = 0$ in a
  distributional sense. Moreover, it does not change sign and is therefore
  the ground state of the Dirichlet Laplacian on $D_i$.

  Extending each $\phi_{*,i}$ by zero outside its domain, we obtain
  $k$ linearly independent eigenfunctions $\widetilde{\phi_{*,i}}$ of
  $-\Delta^\chi_\infty$ corresponding to the eigenvalue $\lambda_*$.
  Conversely, for any $\phi \in \Phi$, its restriction $\phi_i$ is a
  $\lambda_*$-eigenfunction of the Dirichlet Laplacian on $D_i$ (if
  non-zero), and therefore must be proportional to the ground state.
  We conclude that $\dim\Phi = k$.

  From Proposition~\ref{prop:Deltachi} we get $\gN \phi_* = 0$.  Let
  $\psi \in \Phi$ be another function such that $\gN \psi = 0$.  Since
  the restriction of $\psi$ to (say) subdomain $D_1$ is a multiple of
  its ground state, there is a linear combination of $\phi_*$ and
  $\psi$ which identically vanishes on $D_1$.  By a straightforward
  extension of the unique continuation principle to $\Delta^\chi$,
  this linear combination is zero everywhere and therefore $\psi$ is a
  multiple of $\phi_*$. 
  
  Next, Corollary~\ref{cor:equalangle} implies that each $D_i$ has piecewise 
  smooth boundary with convex corners, so it follows from \cite[Remark~3.2.4.6]{Gr} 
  that $\phi_i \in H^2(D_i)$ for any $\phi \in \Phi$.

  Finally, let
  \begin{equation}
    \label{eq:def_extension_by_zero}
    E_i \colon L^2(\Gamma_i) \to L^2(\Gamma)  
  \end{equation}
  denote extension by zero.  We have
  \begin{equation}
    \label{eq:extension_trace}
    \gN \widetilde{\phi_{*,i}} = E_i\left(\chi_i \frac{\p \phi_{*,i}}{\p\nu_i}\right),
  \end{equation}
  and therefore the claim $\gN(\Phi) \subset \HGc$
  follows from the next proposition applied to $\phi_{*,i}$.
\end{proof}

\begin{proposition}
  \label{normalreg}
  If $u \in H^2(D_i) \cap H^1_0(D_i)$, then $\p_{\nu_i}u \in
  H^{1/2}(\Gamma_i)$, $\chi_i \p_{\nu_i}u \in H^{1/2}(\Gamma_i)$ and
  $E_i(\chi_i \p_{\nu_i}u) \in H^{1/2}_\chi(\Gamma)$.
\end{proposition}

The assumption that $u$ vanishes on the boundary is essential. If $D_i$ 
has corners, then the unit normal $\nu_i$ is discontinuous there, and for a 
general function in $H^2(D_i)$, or even $C^\infty(\overline D_i)$, there is no guarantee that
$\p_{\nu_i} u \in H^{1/2}(\p D_i)$. A simple example is $u(x,y) = x$ on the 
unit square; its normal derivative is piecewise constant, but 
is not contained in $H^{1/2}$.

Localizing around a single corner and performing a suitable change of variables, it suffices to prove the result for the model domain $D = \bbR_+ \times \bbR_+$, which has boundary $\p D = \overline{(\bbR_+ \times \{0\}) \cup (\{0\} \times \bbR_+)}$. We first recall some preliminary results on boundary Sobolev spaces.

\begin{lemma}\cite[Theorem~1.5.2.3]{Gr}
\label{lem:glue}
Given $f_1, g_1 \in H^{1/2}(\bbR_+)$, the composite function
\[
	h = \begin{cases} f_1, & \text{ on } \bbR_+ \times \{0\} \\
	g_1, & \text{ on } \{0\} \times \bbR_+ \end{cases}
\]
is in $H^{1/2}\big(\p D)$ if and only if
\begin{equation}
\label{compat0}
	\int_0^1 |f_1(t) - g_1(t)|^2 \frac{dt}{t} < \infty.
\end{equation}
\end{lemma}

In particular, the conclusion $h \in H^{1/2}(\p D)$ holds for any $f_1$ and $g_1$ satisfying the stronger condition
\begin{equation}
\label{compat1}
	\int_0^1 |f_1(t)|^2 \frac{dt}{t} + \int_0^1 |g_1(t)|^2 \frac{dt}{t} < \infty.
\end{equation}
We also need to know the image of the trace map on each smooth component of the boundary.

\begin{lemma}\cite[Theorem~1.5.2.4]{Gr}
\label{lem:trace}
The trace map
\begin{align}
\label{trace}
\begin{split}
	H^2(\bbR_+ \times \bbR_+) &\rightarrow H^{3/2}(\bbR_+) \times H^{1/2}(\bbR_+) \times H^{3/2}(\bbR_+) \times H^{1/2}(\bbR_+) \\
	u(x,y) &\mapsto \left( u(x,0), \frac{\p u}{\p y}(x,0),  u(0,y), \frac{\p u}{\p x}(0,y) \right)
\end{split}
\end{align}
is continuous, with image consisting of all $(f_0,f_1,g_0,g_1)$ that satisfy the compatibility conditions
\begin{equation}
\label{compat2}
	f_0(0) = g_0(0)
\end{equation}
and
\begin{equation}
\label{compat3}
	\int_0^1 \frac{|f_0'(t) - g_1(t)|^2}{t} dt + \int_0^1 \frac{|f_1(t) - g_0'(t)|^2}{t} dt < \infty.
\end{equation}
\end{lemma}

While the above two lemmas are both if and only if statements, we do not require their full strength in the following proof. It is enough to know that~\eqref{compat1} is a \emph{sufficient} condition for $h$ to be in $H^{1/2}(\p D)$, and~\eqref{compat3} is a \emph{necessary} condition for $(f_0,f_1,g_0,g_1)$ to be in the image of the trace map defined in~\eqref{trace}.

\begin{proof}[Proof of Proposition~\ref{normalreg}]
As mentioned above, it suffices to prove the result for the model domain $D = \bbR_+ \times \bbR_+$. If $u \in H^2(D) \cap H^1_0(D)$, then its corresponding traces
\[
	\big(f_0(x), f_1(x), g_0(y), g_1(y)\big) = \left( u(x,0), \frac{\p u}{\p y}(x,0),  u(0,y), \frac{\p u}{\p x}(0,y) \right)
\]
satisfy $f_0(x) = 0$ and $g_0(y) = 0$ for all $x,y > 0$, so Lemma~\ref{lem:trace} implies
\[
	\int_0^1 \left| \frac{\p u}{\p x}(0,t)\right|^2 \frac{dt}{t} + \int_0^1 \left| \frac{\p u}{\p y}(t,0)\right|^2 \frac{dt}{t} < \infty.
\]
It then follows from Lemma~\ref{lem:glue} that the function
\[
	h
	= \begin{cases} \frac{\p u}{\p y}(x,0), & \text{ on } \bbR_+ \times \{0\} \\
	 \frac{\p u}{\p x}(0,y), & \text{ on } \{0\} \times \bbR_+ \end{cases}
\]
is in $H^{1/2}(\p D)$. Since $h = -\p_\nu u$, this proves the first part of the proposition.

Since the weight $\chi$ is constant (either $+1$ or $-1$) on each axis, we obtain
\begin{equation}
\label{finite}
	\int_0^1 \left| \chi \frac{\p u}{\p x}(0,t)\right|^2 \frac{dt}{t} = \int_0^1 \left| \frac{\p u}{\p x}(0,t)\right|^2 \frac{dt}{t} < \infty,
\end{equation}
and similarly for the integral involving $\p u/\p y$, which implies $\chi \p_\nu u \in H^{1/2}(\p D)$.

This completes our analysis on the model domain $D$, where we have shown that $\chi_i \p_{\nu_i}u \in H^{1/2}(\Gamma_i)$. Finally, we consider the extension $E_i(\chi_i \p_{\nu_i}u)$ to the rest of $\Gamma$. Fixing another domain $D_j$, we must show that $\chi_j E_i(\chi_i \p_{\nu_i}u) \in H^{1/2}(\Gamma_j)$. On each smooth segment of $\Gamma_j$, this function is given by $\chi_j \chi_i \p_{\nu_i}u$ (if $\Gamma_i$ intersects $\Gamma_j$ nontrivially) and $0$ otherwise. Either way, it follows from~\eqref{finite} that the finiteness condition~\eqref{compat1} holds, and so $\chi_j E_i(\chi_i \p_{\nu_i}u) \in H^{1/2}(\Gamma_j)$, as was to be shown.
\end{proof}


\section{Defining the weighted Dirichlet-to-Neumann operator}
\label{sec:def}
In this section we construct the weighted, two-sided Dirichlet-to-Neumann operator $\DN(\Gamma,\lambda_*,\chi)$ for a $\chi$-nodal partition with eigenvalue $\lambda_*$. In Section~\ref{sec:formdef} we give a definition using the standard theory of self-adjoint operators and coercive bilinear forms; the details of this construction are then provided in Sections~\ref{sec:Sprelim} and~\ref{sec:formproof}.

As mentioned in the introduction, the construction is rather involved because $\lambda_*$ is a Dirichlet eigenvalue on each $D_i$. In this case one can also view the Dirichlet-to-Neumann map as a multi-valued operator (or linear relation); this approach is described in \cite{AEKS,AM,BtE}. Another difficulty is that $\Gamma$ has corners. While the Dirichlet-to-Neumann map can be defined on domains with minimal boundary regularity (see \cite{AE11}), our results require delicate regularity properties, as in Proposition~\ref{prop:regDirichlet}, that are not available in that case.

\subsection{Definition via bilinear forms}
\label{sec:formdef}
We define $\DN(\Gamma,\lambda_*,\chi)$ as the self-adjoint operator corresponding to a bilinear form on the closed subspace
\begin{equation}
\label{Sdef}
	S_\chi := \left\{g\in L^2(\Gamma) : \int_{\Gamma_i} \chi_i g_i  \frac{\p \phi_{*,i}}{\partial \nu_i } = 0, \ i=1,\dots, k\right\}
\end{equation}
of $L^2(\Gamma)$, where $\phi_{*,i}$ denotes the restriction of $\phi_{*}$ to $D_i$, and we recall that $\Gamma_i = \overline{\p D_i \cap \Omega}$.

Let $g \in \HGc \cap S_\chi$. For each $i$, the problem
\begin{equation}\label{eq:Di2}
	\begin{cases}
	-\Delta u_i = \lambda_* u_i &\text{in $D_i$},\\
		u_i = \chi_i g_i & \text{on } \Gamma_i,\\
		u_i=0 & \text{on } \partial D_i\setminus \Gamma_i,
	\end{cases}
\end{equation}
has a unique solution $u_i^g \in H^1(D_i)$ that satisfies the orthogonality condition
\begin{equation}\label{eq:can}
	\int_{D_i}u_i^g \,\phi_{*,i} = 0.  
\end{equation}
Using these solutions, we define the symmetric bilinear form
\begin{equation}
\label{adef}
	a(g,h) := \sum_i \int_{D_i} \big(\nabla u_i^g \cdot \nabla u_i^h - \lambda_* u_i^g u_i^h \big), \qquad \dom(a) = \HGc \cap S_\chi.
\end{equation}
It follows from Lemma~\ref{lemma:Ntrace} that the two-sided normal derivative $\gN u^g \in \HmGc$ is defined.

The main result of this section is the following.

\begin{theorem}
\label{thm:Aform}
Let $\cD = \{D_i\}$ be a $\chi$-nodal partition.
\begin{enumerate}
\item The bilinear form  $a$ defined in~\eqref{adef} generates a self-adjoint operator
  $\DN(\Gamma,\lambda_*,\chi)$, which has domain
  \begin{equation}
    \label{Adom}
    \dom\big(\DN(\Gamma,\lambda_*,\chi)\big)
    = \big\{g \in \HGc \cap S_\chi \colon \gN u^g \in L^2(\Gamma) \big\},
  \end{equation}
  and is given by
  \begin{equation}
    \label{Aexplicit}
    \DN(\Gamma,\lambda_*,\chi)g = \Pi_\chi(\gN u^g),
  \end{equation}
  where $\Pi_\chi$ is the $L^2(\Gamma)$-orthogonal projection onto
  $S_\chi$.
\item For each $g \in \dom(\DN(\Gamma,\lambda_*,\chi))$, there exists a function
  $\tilde u \in \dom(t^\chi)$ such that $\gN \tilde u \in S_\chi$ and
  $\tilde u_i$ solves~\eqref{eq:Di2} for each $i$, hence
  \begin{equation}
    \label{Aexplicit2}
    \DN(\Gamma,\lambda_*,\chi)g = \gN \tilde u.
  \end{equation}
  If we additionally require $\int_\Omega \tilde u \phi_* = 0$, then $\tilde u$ is unique.

\item For any other set of valid weights $\{\tilde\chi_i\}$, we have:
  \begin{enumerate}
  \item If $\chi$ and $\tilde\chi$ are edge equivalent, then
    $\DN(\Gamma,\lambda_*,\chi)$ is unitarily equivalent to
    $\DN(\Gamma,\lambda_*,\tilde\chi)$;
  \item If $\chi$ and $\tilde\chi$ are domain equivalent, then
    $\DN(\Gamma,\lambda_*,\chi) = \DN(\Gamma,\lambda_*,\tilde\chi)$.
  \end{enumerate}
  Consequently, $\DN(\Gamma,\lambda_*,\chi)$ and
  $\DN(\Gamma,\lambda_*,\tilde\chi)$ are unitarily equivalent for any
  two valid sets of weights $\chi$ and $\tilde\chi$.
\end{enumerate}
\end{theorem}

\begin{remark}
\label{rem:ug}
If $u \in \dom(t^\chi)$ and $u_i$ solves~\eqref{eq:Di2} for each $i$, it must be of the form $u = u^g + \phi$ for some $\phi \in \Phi$, by Proposition~\ref{prop:regDirichlet}. Since $\gN\phi \in S_\chi^\perp$, we have $\Pi_\chi (\gN u^g) = \Pi_\chi (\gN u)$, meaning $u^g$ can be replaced by any other solution to~\eqref{eq:Di2}. 
The distinguished solution $u_i^g$ has nice analytic properties, which we will use in Lemma~\ref{lemma:bounds} to prove that $a$ is semi-bounded and closed. On the other hand, $\gN u^g$ may not be in the subspace $S_\chi$, so we need to apply the orthogonal projection $\Pi_\chi$ in~\eqref{Aexplicit}. By choosing different solutions to~\eqref{eq:Di2} we can eliminate this projection, as in~\eqref{Aexplicit2}.
\end{remark}

\subsection{The subspace $S_\chi$}
\label{sec:Sprelim}
We start by discussing some useful properties of the subspace $S_\chi$ defined in~\eqref{Sdef}.
Recall that $\Phi$ is the kernel of the Dirichlet Laplacian $\Delta^\chi_\infty$, as described in Proposition~\ref{prop:regDirichlet}.

\begin{lemma}
  \label{S:explicit}
  The
  subspace $S_\chi \subset L^2(\Gamma)$ can be written as
  \begin{equation}
    \label{Schi}
    S_\chi = \left\{g\in L^2(\Gamma) \colon
      \left<\gN \phi, g \right>_{L^2(\Gamma)} = 0 \text{ for all } \phi \in
      \Phi\right\}
    = \big( \gN(\Phi) \big)^\perp.
  \end{equation}
Therefore, it is a closed subspace of codimension $k-1$.
\end{lemma}

\begin{proof}
  Formula~\eqref{Schi} is a direct consequence of the properties of
  the space $\Phi$ and equation~\eqref{eq:extension_trace}, since
  \[
    \int_{\Gamma_i} \chi_i g_i  \frac{\p \phi_{*,i}}{\p\nu_i} =
    \int_\Gamma g \gN \widetilde{\phi_{*,i}}.
  \]
  From Proposition~\ref{prop:regDirichlet} we have
  \[
	\dim \gN(\Phi) = \dim \Phi - \dim\ker \big(\gN\big|_\Phi\big) = k-1.
  \]
  Since $S_\chi^\perp = \gN(\Phi)$, this completes the proof.
\end{proof}

\begin{lemma}
\label{lemma:dense}
The set $\HGc \cap S_\chi$ is dense in $S_\chi$.
\end{lemma}

\begin{proof}
We first claim that $\HGc$ is dense in $L^2(\Gamma)$. Fix $g \in L^2(\Gamma)$ and let $\epsilon>0$. Letting $\tilde\Gamma$ denote the smooth part of $\Gamma$, which is diffeomorphic to a finite number of open intervals, we can find a function $\tilde g$ on $\Gamma$ such that $\|g - \tilde g\|_{L^2(\Gamma)} < \epsilon$ and $\tilde g \in C^\infty_0(\tilde\Gamma)$. Since the weights $\chi_i$ are constant on each component of $\tilde\Gamma$, it follows that $\chi_i \tilde g \in C^\infty_0(\tilde\Gamma \cap \Gamma_i)$, and hence $\chi_i \tilde g \in H^{1/2}(\Gamma_i)$,  for each $i$. This implies $\tilde g \in \HGc$ and thus proves the claim.

Now let $g \in S_\chi$, $\epsilon>0$, and choose $\tilde g \in \HGc$
as above.  Lemma~\ref{S:explicit} implies
$(I - \Pi_\chi)\tilde g \in \gN(\Phi)$, which in turn belongs to $\HGc$
by Proposition~\ref{normalreg}.  Therefore,
\begin{equation}
  \label{eq:proj_in_Hhalf}
  \Pi_\chi \tilde g = \tilde g - (I - \Pi_\chi)\tilde g \in \HGc \cap S_\chi. 
\end{equation}
We now use the fact that $\Pi_\chi g = g$ to obtain
\[
	\| g - \Pi_\chi \tilde g \|_{S_\chi} = \| \Pi_\chi(g - \tilde g) \|_{S_\chi} \leq 
	\| g - \tilde g \|_{L^2(\Gamma)} < \epsilon,
\]
as was to be shown.
\end{proof}

Finally, we describe the set of functionals in $H_\chi^{-1/2}(\Gamma)
= \HGc^*$ that vanish on $\HGc \cap S_\chi$. This will be
used below, in the proof of Theorem~\ref{thm:Aform}, when we describe
the domain of the Dirichlet-to-Neumann map.

\begin{lemma}
  \label{lemma:vanish}
  If $\tau \in H^{-1/2}_\chi(\Gamma)$ and $\tau(g) = 0$ for all $g \in \HGc \cap S_\chi$, then there exists a function $h \in \gN(\Phi) = S_\chi^\perp$ such that $\tau(g) = \left<g,h\right>_{L^2(\Gamma)}$  for all $g \in \HGc$.

\end{lemma}

\begin{proof}
  From Lemma~\ref{lemma:dense} we have the $L^2(\Gamma)$-orthogonal
  decomposition
  \begin{equation*}
    \HGc = \left(\HGc \cap S_\chi\right) \oplus \gN(\Phi),
  \end{equation*}
  therefore any
  functional that vanishes on $\HGc \cap S_\chi$ is a functional on $\gN(\Phi)$
  extended by zero. Since $\gN(\Phi)$ is finite dimensional, a functional 
  $\hat\tau \colon \gN(\Phi) \to \bbR$ is continuous with respect to any choice of 
  norm. In particular, it is continuous with respect to the $L^2(\Gamma)$ 
  norm, so there exists $h \in \gN(\Phi)$ such that $\hat\tau(g) =  
  \left<g,h\right>_{L^2(\Gamma)}$  for all $g \in \gN(\Phi)$.
\end{proof}

\subsection{Proof of Theorem~\ref{thm:Aform}}
\label{sec:formproof}
From Lemma~\ref{lemma:dense} we know that the symmetric bilinear form
$a$ is densely defined. The next step is to show that it is
semi-bounded and closed. This is an immediate consequence of the 
completeness of $\HGc$ and the
following inequalities; see, for instance, \cite[Section~11.2]{Schm12}. 

\begin{lemma}
\label{lemma:bounds}
There exist constants $C, c > 0$ and $m \in \bbR$ such that
\begin{equation}
\label{a:bound2}
	|a(g,h)| \leq C \|g\|_{\HGc} \|h\|_{\HGc}
\end{equation}
and
\begin{equation}
\label{a:coercive2}
	a(g,g) \geq c \|g\|^2_{\HGc} + m\| g\|^2_{L^2(\Gamma)}
\end{equation}
for all $g,h \in \HGc \cap S_\chi$.
\end{lemma}

In the proof we let $C,c$ denote positive constants, and $m$ a real constant, whose meaning may change from line to line.

\begin{proof}
For each $i$ the unique solution $u_i^g$ to~\eqref{eq:Di2} and~\eqref{eq:can} satisfies a uniform estimate
\[
  \|u_i^g\|_{H^1(D_i)} \leq  C \big\| \cE_i(\chi_i g_i) \big\|_{H^{1/2}(\p D_i)}.
\]
Recalling the definition of the $\HGc$ norm in~\eqref{norm:HGc}, it follows that
\[
	|a(g,h)| \leq C \|g\|_{\HGc} \|h\|_{\HGc}
\]
for all $g,h \in \HGc \cap S_\chi$.

On the other hand, a standard compactness argument
(see~\cite[Lemma~2.3]{AM12}) shows that for any $\epsilon>0$ there exists a constant
$K(\epsilon) > 0$ such that
\begin{equation}
  \label{ubound}
  \|u_i\|_{L^2(D_i)}^2 \leq \epsilon \| \nabla u_i \|_{L^2(D_i)}^2 + K(\epsilon) \big\| u_i|_{\Gamma_i} \big\|_{L^2(\Gamma_i)}^2
\end{equation}  
for all $u_i$ in the set
\[
	\left\{u_i \in H^1(D_i) : \Delta u_i + \lambda_* u_i = 0, \ \int_{D_i} u_i \phi_{*,i} = 0, \ u_i\big|_{\p D_i \cap \p\Omega} = 0 \right\}.
\]
In particular, the estimate~\eqref{ubound} holds for each $u_i^g$. It
then follows, exactly as in~\cite[Proposition~3.3]{AM12}, that 
\begin{align*}
	\int_{D_i} \big(|\nabla u_i^g|^2 - \lambda_* |u_i^g|^2 \big) &\geq \frac12 \|u_i^g\|_{H^1(D_i)}^2 + m \| g_i \|^2_{L^2(\Gamma_i)}\\
	& \geq c \big\| \cE_i(\chi_i g_i) \big\|^2_{H^{1/2}(\p D_i)} + m \| g_i \|^2_{L^2(\Gamma_i)}
\end{align*}
for each $i$, with constants $c > 0$ and $m \in \bbR$, and hence
\[
	a(g,g) \geq c \|g\|^2_{\HGc} + m \| g\|^2_{L^2(\Gamma)}
\]
for all $g \in \HGc \cap S_\chi$.
\end{proof}

We are now ready to prove the main result.

\begin{proof}[Proof of Theorem~\ref{thm:Aform}]  
From Lemmas~\ref{lemma:dense} and~\ref{lemma:bounds} we know that the
symmetric bilinear form $a$ is densely defined, lower semi-bounded and
closed, so it generates a self-adjoint operator on $S_\chi$, which we
denote by $A$ for brevity. Its domain is given by
\begin{align}
\begin{split}\label{Adom1}
	\dom(A) = \big\{ g \in & \HGc \cap S_\chi : \text{there exists } f \in S_\chi \text{ such that } \\ &a(g,h) = \left<f,h\right>_{L^2(\Gamma)} \text{ for all } h \in \HGc \cap S_\chi \big\},
\end{split}
\end{align}
and $Ag = f$ for any such $g$.

We now characterize the domain of $A$. First suppose $g \in \dom(A)$,
and let $f = Ag \in S_\chi$. Using Lemma~\ref{lemma:Ntrace} and the
definition of $a$ in~\eqref{adef}, we get
\begin{equation}
\label{green}
	a(g,h) = \left<\gN u^g, h \right>_*
\end{equation}
for all $h \in \HGc \cap S_\chi$. On the other hand,~\eqref{Adom1} implies
\[
	a(g,h) = \left<f,h\right>_{L^2(\Gamma)} = \int_\Gamma fh,
\]
so we find that
\[
	\gN u^g - f \in H_\chi^{-1/2}(\Gamma)
\]
vanishes on $\HGc \cap S_\chi$. From Lemma~\ref{lemma:vanish} we get
$\gN u^g - f \in \gN(\Phi) = S_\chi^\perp$, and hence $\gN u^g \in L^2(\Gamma)$. Since $f \in S_\chi$, it
follows that $f = \Pi_\chi(\gN u^g)$.

Conversely, if $g \in \HGc \cap S_\chi$ and $\gN u^g \in L^2(\Gamma)$, we have
\[
	\int_\Gamma \big(\Pi_\chi(\gN u^g)\big) h =  \int_\Gamma (\gN u^g) h = a(g,h)
\]
for all $h \in \HGc \cap S_\chi$. According to~\eqref{Adom1}, this implies $g \in \dom(A)$, with $Ag = \Pi_\chi(\gN u^g)$.

Next, we prove the existence of $\tilde u$. Since $\Pi_\chi$ is the orthogonal projection onto $S_\chi = \gN(\Phi)^\perp$, we have
\[
	Ag = \Pi_\chi(\gN u^g) = \gN u^g - \gN\phi
\]
for some $\phi\in\Phi$.  Setting $\tilde u = u^g - \phi$, we
obtain $Ag = \gN \tilde u$, as required.  If $\hat u$ is another function in $\dom(t^\chi)$ such that
$Ag = \gN \hat u$ and $\hat u_i$ solves~\eqref{eq:Di2} for each
$i$, then $\hat u - \tilde u \in \Phi$ and also
$\hat u - \tilde u \in \ker \gN$.  By
Proposition~\ref{prop:regDirichlet}, $\hat u - \tilde u$ is a multiple
of $\phi_*$, and so requiring $\tilde u$ to be orthogonal to $\phi_*$
determines it uniquely.

Finally, we establish the dependence on the weights. If $\chi$ and $\tilde\chi$ are edge equivalent, the desired unitary transformation is multiplication by $\tilde\chi_i / \chi_i$ on $\Gamma_i$. The edge equivalence ensures this is well-defined, since $\tilde\chi_i / \chi_i = \tilde\chi_j / \chi_j$ on $\Gamma_i \cap \Gamma_j$. The result when $\chi$ and $\tilde\chi$ are domain equivalent follows immediately from the definition.
\end{proof}

\section{The spectral flow: proof of Theorem~\ref{thm:general}}
\label{sec:flow}

To prove our main theorem we study the spectral flow of a family of self-adjoint operators. This idea was pioneered by Friedlander in \cite{F91}, though our approach is closer to that of~\cite{AM12,AM}. To characterize the negative eigenvalues of $\DN(\Gamma,\lambda_*,\chi)$ it is fruitful to study the family of operators $-\Delta^\chi_\sigma$, $0\leq\sigma<\infty$, induced by the symmetric bilinear form
\begin{equation}
\label{form:tsigma}
	t^\chi_\sigma(u,v)=	
	\sum_{i=1}^k \int_{D_i} \nabla u_i \cdot \nabla v_i + \sigma\int_\Gamma uv,
	\qquad \dom(t^\chi_\sigma) = \dom(t^\chi),
\end{equation}
where $\dom(t^\chi)$ was defined
in~\eqref{eq:domt_H1}--\eqref{eq:domt_Dint}.  As in
Proposition~\ref{prop:Deltachi}, it can be shown that each
$\Delta^\chi_\sigma$ is self-adjoint, with domain
\begin{align}
\label{domDsigma}
\begin{split}
	\dom(\Delta^\chi_\sigma) = \big\{ u \in \dom(t^\chi) : \Delta u_i \in L^2(D_i) \text{ for each } i \text{ and } \gN u + \sigma\gD u = 0  \big\}.
\end{split}
\end{align}
It can be easily seen that the eigenfunction $\phi_*$ of
$\Delta^{\chi}$ that vanishes on the set $\Gamma$ is an eigenfunction
of $\Delta^\chi_\sigma$ for all $\sigma$.  We can therefore consider the
reduced operator $\widehat\Delta^{\chi}_\sigma$, which is simply
$\Delta^\chi_\sigma$ restricted to $\spn\{\phi_*\}^\perp$. We recall
(see Section~\ref{sec:regularity}) that $\Delta^\chi_\infty$ is the
Laplacian on $\Omega$ with Dirichlet boundary conditions imposed on
$\partial\Omega\cup\Gamma$.

\begin{proposition}\label{prop:DtoNeig}
For each $\sigma\in[0,\infty)$ the linear mapping
\begin{equation}
\label{linear}
\begin{gathered}
	T \colon \ker(\Delta^\chi_\sigma + \lambda_*) \longrightarrow \ker \big(\DN(\Gamma,\lambda_*,\chi) + \sigma \big), \qquad Tu = \gD u,
\end{gathered}
\end{equation}
is surjective, and its kernel is spanned by $\phi_*$, hence
\[
	\dim\ker(\Delta^\chi_\sigma + \lambda_*) - 1 = \dim\ker\big(\DN(\Gamma,\lambda_*,\chi) + \sigma\big).
\]
Equivalently, in terms of the reduced operator, the restriction of $T$ to $ \ker(\widehat\Delta^{\chi}_\sigma + \lambda_*)$ is bijective and
\[
	\dim\ker(\widehat\Delta^{\chi}_\sigma + \lambda_*) = \dim\ker\big(\DN(\Gamma,\lambda_*,\chi)+\sigma\big).
\]
\end{proposition}

\begin{proof}
We first show that $T$ is well-defined. Assume that $u$ is an eigenfunction of $-\Delta^\chi_\sigma$ with eigenvalue $\lambda_*$. From~\eqref{domDsigma} we see that $u$ satisfies the transmission condition $\gN u + \sigma \gD u = 0$ on $\Gamma$. On each $D_i$ we can use Green's second identity to conclude that
\[
	0 = \int_{\Gamma_i}u_i \frac{\p \phi_{*,i}}{\partial \nu_i} = \int_{\Gamma_i} \chi_i (\chi_i u_i) \frac{\p \phi_{*,i}}{\partial \nu_i}.
\]
This means that the Dirichlet trace $\gD u \in \HGc$ belongs to the subspace $S_\chi$ defined in~\eqref{Sdef}. Moreover, since $\gN u = -\sigma \gD u$ is contained in $L^2(\Gamma)$, we see from~\eqref{Adom} that $\gD u$ belongs to the domain of $\DN(\Gamma,\lambda_*,\chi)$, with
\begin{equation}
\label{DNproj}
	\DN(\Gamma,\lambda_*,\chi) \gD u = \Pi_\chi(\gN u) = -\sigma \Pi_\chi(\gD u) = -\sigma\gD u.
\end{equation}
This means $\gD u \in \ker\big(\DN(\Gamma,\lambda_*,\chi)+\sigma\big)$, so $T$ is well-defined.

We next show that $T$ is surjective. Let $g \in \ker\big(\DN(\Gamma,\lambda_*,\chi)+\sigma\big)$ be given. From the second part of Theorem~\ref{thm:Aform}, we know that there exists $\tilde u_i \in H^1(D_i)$ satisfying the equation $\Delta \tilde u_i + \lambda_* \tilde u_i = 0$ and the boundary conditions $\gD \tilde u = g$, such that
\begin{equation}
	\gN \tilde u = \DN(\Gamma,\lambda_*,\chi)g = -\sigma g.
\end{equation}
This is precisely the transmission condition $\gN \tilde u + \sigma \gD \tilde u = 0$, so we conclude from~\eqref{domDsigma} that $\tilde u \in \dom(\Delta_\sigma^\chi)$ and hence $\tilde u \in \ker(\Delta_\sigma^\chi + \lambda_*)$. Since $T\tilde u = \gD \tilde u = g$, this proves surjectivity.

It remains to prove that the kernel of $T$ is spanned by $\phi_*$. From Proposition~\ref{prop:regDirichlet} we know that $\phi_*\in \ker(\Delta_\sigma^\chi + \lambda_*)$ for all $\sigma$ and $\gD \phi_* = 0$, hence $\phi_* \in \ker T$. Finally, suppose that $u$ is any function in $\ker T$. This means $u \in \ker(\Delta_\sigma^\chi + \lambda_*)$ and $\gD u = 0$, hence $\gN u = 0$ by the transmission condition, so it follows from Proposition~\ref{prop:regDirichlet} that $u$ is proportional to $\phi_*$.
\end{proof}

\begin{remark}
In the above proof, in particular~\eqref{DNproj}, we see that if $u \in \ker(\Delta_\sigma^\chi + \lambda_*)$, so that $\gD u$ is an eigenfunction of $\DN(\Gamma,\lambda_*,\chi)$, then $\gN u \in S_\chi$, and hence $\Pi_\chi(\gN u) = \gN u$. In other words, it is the particular solution $\tilde u$ whose existence is guaranteed by the second part of Theorem~\ref{thm:Aform}.
\end{remark}

We are now ready to prove our main result.

\begin{proof}[Proof of Theorem~\ref{thm:general}]
The equality~\eqref{eq:mult2} follows from Proposition~\ref{prop:DtoNeig} with $\sigma=0$. To prove~\eqref{eq:def} we consider the spectral flow for the reduced operator family $-\widehat\Delta^\chi_\sigma$ defined above. Since this is an analytic family of self-adjoint operators for $0 \leq \sigma < \infty$, we can arrange the eigenvalues into analytic branches $\{\gamma_m(\sigma)\}$ such that:
\begin{enumerate}
	\item $\{\gamma_m(0)\}$ are the ordered eigenvalues of $-\widehat\Delta^\chi_0$, repeated according to multiplicity;
	\item each function $\sigma\mapsto \gamma_m(\sigma)$ is non-decreasing;
	\item as $\sigma \to \infty$, the $\gamma_m(\sigma)$ converge to the eigenvalues of $-\widehat\Delta^\chi_\infty$.
\end{enumerate}
The first statement is simply our convention for labelling the branches, the second follows from the 
monotonicity of the quadratic form $t^\chi_\sigma(u,u)$ from \eqref{form:tsigma}, and the third can be proved using the method of \cite[Theorem~2.5]{AM12}.

At $\sigma=0$ the operator $\widehat\Delta^\chi_0$ has $\ell-1$ eigenvalues below $\lambda_*$. On the other hand, at $\sigma=\infty$ the first eigenvalue of $-\Delta^\chi_\infty$ is $\lambda_*$, with multiplicity $k$ (one for each nodal domain). This means the first eigenvalue of the reduced operator $-\widehat\Delta^\chi_\infty$ is also $\lambda_*$, but with multiplicity $k-1$.

Therefore, of the first $\ell-1$ eigenvalue curves, precisely $k-1$ converge to $\lambda_*$, while the remaining $\ell-k$ converges to strictly larger values, and hence intersect $\lambda_*$ at some finite value of $\sigma$. In other words,
\[
	\ell - k = \# \{m : \gamma_m(\sigma) = \lambda_* \text{ for some } \sigma \in (0,\infty) \}.
\]
From Proposition~\ref{prop:DtoNeig} we know that $\lambda_*$ is an eigenvalue of $-\widehat\Delta^\chi_\sigma$ if and only if $-\sigma$ is an eigenvalue of $\DN(\Gamma,\lambda_*,\chi)$, with the same multiplicity, and hence
\[
	\# \{m : \gamma_m(\sigma) = \lambda_* \text{ for some } \sigma \in (0,\infty) \}
	= \Mor \DN(\Gamma,\lambda_*,\chi)
\]
is the number of negative eigenvalues of $\DN(\Gamma,\lambda_*,\chi)$, counted with multiplicity.
\end{proof}

\section{Equipartitions of the unit circle revisited}
\label{sec:circle}
Here we analyze equipartitions of the circle, calculating explicitly the different terms in Theorem~\ref{thm:general}. The same example was previously considered in~\cite{HeSu}, but with the Dirichlet-to-Neumann map evaluated at $\lambda_* + \epsilon$, as described in the introduction. Also, in~\cite{HeSu} the magnetic point of view was used, with the operator $T=-\bigl(\frac d {d\theta} -\frac{\i}{2}\bigr)^2$. We use here the equivalent presentation with cuts, replacing $T$ by $-\Delta^\chi$.

Leting $\mathcal D=\{D_i\}_{i=1}^k$ be a $k$-equipartition of the circle, we will show that $\cD$ is $\chi$-nodal, corresponding to a $\Delta^\chi$ eigenvalue of multiplicity two, with defect $\delta(\cD) = 0$. Comparing with Theorem~\ref{thm:general}, we should thus have
\begin{equation}
\label{circle}
	\Mor  \DN(\Gamma,\lambda_*,\chi) = 0, \qquad
	\dim\ker \DN(\Gamma,\lambda_*,\chi) = 1.
\end{equation}
Indeed, we find that $\DN(\Gamma,\lambda_*,\chi)$ is identically zero on the space $S_\chi$, which is one dimensional, confirming \eqref{circle}.

\begin{remark}
Recall that $S_\chi \subset L^2(\Gamma)$ has codimension $k-1$. A $k$-partition of the circle has $k$ boundary points, so $L^2(\Gamma) \cong \bbR^k$ and hence $S_\chi$ is one dimensional. A $k$ partition of an interval, however, has only $k-1$ boundary points, and so $S_\chi$ is zero dimensional. In this case the nullity and Morse index of $\DN$ must be zero, so Theorem~\ref{thm:general} says that the partition has zero deficiency and corresponds to a simple eigenvalue, thus reproducing the Sturm oscillation theorem.
\end{remark}

We view the circle as $[0,2\pi]$ with the endpoints identified.
We choose as division points $\theta_i=2\pi i/k$ for $0\leq i\leq k$,
naturally identifying $\theta_0$ and $\theta_k$.  The partition
thus consists of the subintervals
\[
  D_i=\bigl(\theta_{i-1},\theta_i\bigr),\quad 1\leq i\leq k,
\]
and the boundary set is given by $\Gamma=\{\theta_i\}_{i=0}^{k-1}$. We
next define the weight functions $\chi_i$, which in our case are
functions on $\partial D_i=\{\theta_{i-1},\theta_i\}$ with values in
$\{\pm 1\}$.

\begin{figure}
  \centering
  \includegraphics{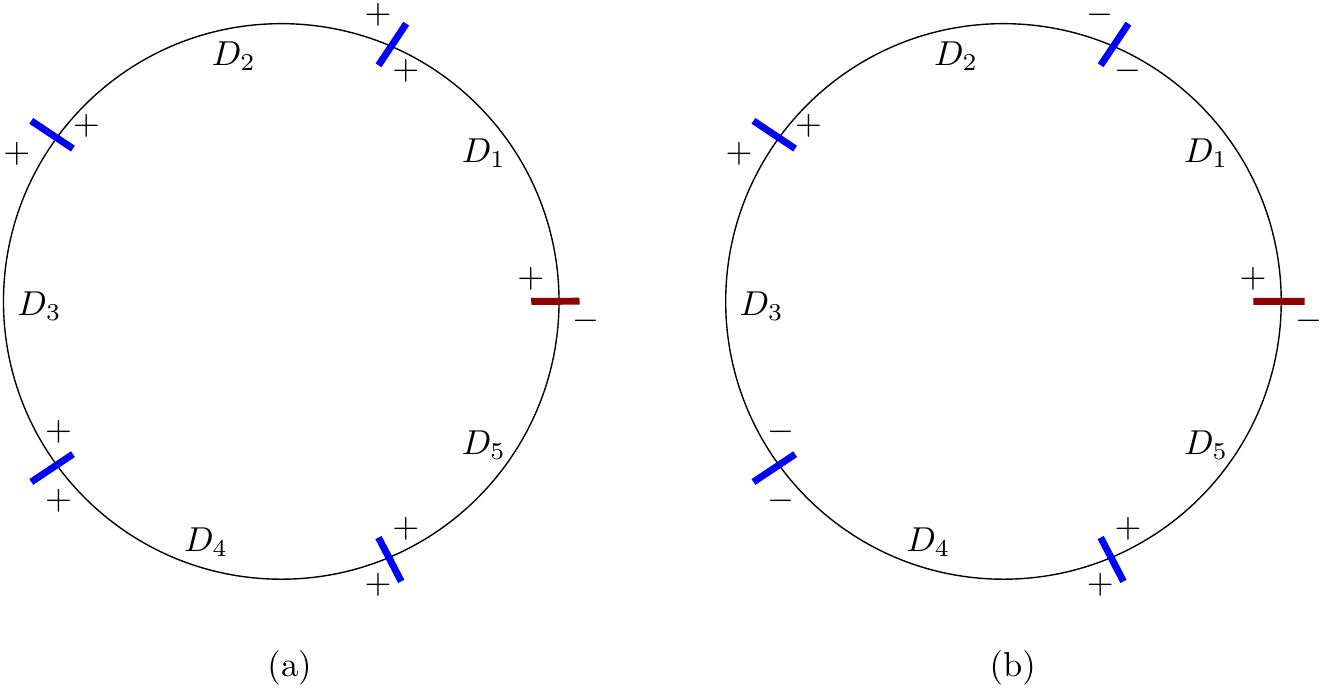}
  \caption{Two choices of valid $\chi$ for an odd partition of a
    circle, as in Section~\ref{sec:circle}.  The value of $\chi_i$ is
    indicated by a $+$ or $-$ next to the corresponding side of the
    partition boundary.  The choices two choices are
    edge-equivalent (Definition~\ref{def:weight_equivalences}).
    One way to see they are valid is to note that they define the same
    ``valid cut'' of $\Gamma$, which is highlighted in red; see Definition~\ref{def:valid_cut} and Proposition~\ref{prop:valid_cuts_valid_weights}.
}
  \label{fig:circle_chi}
\end{figure}

If $k$ is even we are in the bipartite case, and we can choose
$\chi_i \equiv 1$ for each $i$, in which case $\Delta^\chi$ is the
Laplacian. We therefore only consider odd $k$, and introduce a single
cut at $\theta=0$, as was done in~\cite{HeSu}.  As weight functions we
take $\chi_i\equiv 1$ for $0\leq i\leq k-1$, and for $i=k$ we
take $\chi_k(\theta_{k-1})=1$ and $\chi_k(\theta_k)=-1$, see
Figure~\ref{fig:circle_chi}(a).

For the operator $-\Delta^\chi$ we recall from~\eqref{eq:domt_Dint} the compatibility condition $\chi_i u_i=\chi_j u_j$ on the common boundaries of $D_i$ and $D_j$, which here says that functions $u$ in the domain of $-\Delta^\chi$ should be continuous at each~$\theta_i$ except the cut, where $u(0)=-u(2\pi)$. We recall also from~\eqref{domD} the transmission conditions $\chi_i \p_{\nu_i} u_i + \chi_j \p_{\nu_j} u_j = 0$.
The outward normal derivative $\p_{\nu_i}$ is $-\p_\theta$ at the left end-point and $\p_\theta$ at the right end-point, so functions in the domain of $-\Delta^\chi$ should be differentiable at each $\theta_i$ except the cut, where $u'(0)=-u'(2\pi)$. In summary, we have
\begin{equation}
	\dom(\Delta^\chi) = \{u \in H^2(0,2\pi) : u(0) = - u(2\pi), \ u'(0) = -u'(2\pi) \big\}.
\end{equation}
This operator is known as the \emph{anti-periodic Hill operator} or the
magnetic Laplace operator on a circle with flux $1/2$.

The spectrum of~$-\Delta^\chi$ consists of eigenvalues
$\lambda = (j/2)^2$, where $j$ is positive and odd. Each eigenspace is
two dimensional, spanned by $\sin(j\theta/2)$ and $\cos(j\theta/2)$.
The partition $\cD$ is $\chi$-nodal since it is generated by the
eigenfunction $\phi_*(\theta) = \sin(k\theta/2)$.  The minimal label
of the corresponding eigenvalue $\lambda_*=(k/2)^2$ is
$\ell(\mathcal D)=k$ and thus $\delta(\cD) = 0$, as claimed above.

We now turn to the Dirichlet-to-Neumann operator, for which we use a
different valid choice of weights\footnote{The two choices of
  weights are edge equivalent, therefore the Laplacian
  $\Delta^{\hat\chi}$ is identical to $\Delta^\chi$.}, letting 
$\hat\chi_i(\theta_{i-1}) = \cos(k\theta_{i-1}/2) = (-1)^{i-1}$ and
$\hat\chi_i(\theta_{i}) = \cos(k\theta_{i}/2) = (-1)^i $ for all
$1 \leq i \leq k$; see Figure \ref{fig:circle_chi}(b).  The condition for the boundary
data $g = (g_0,g_1,\ldots,g_{k-1})\in\mathbb R^k$ to be in the
subspace $S_{\hat\chi}$ defined in \eqref{Sdef} is $g_{i-1}-g_i = 0$,
yielding $g = h(1,1,\ldots,1)^t$, $h \in \mathbb R$.  For this choice
of $\{\hat\chi_i\}$, the boundary value problem \eqref{eq:Di2} becomes
\begin{equation}
  \label{circleBVP}
  -u_i''=\lambda_* u_i  \ \text{ in $D_i$},
  \quad u_i(\theta_{i-1}) = h \cos(k\theta_{i-1}/2),
  \quad u_i(\theta_i) =  h \cos(k\theta_{i}/2),
\end{equation}
with the general solution
\[
  u_i(\theta) = h \cos(k\theta/2) + c_i\sin(k\theta/2),
\]
where $c_i$ is an arbitrary constant.  According to
Remark~\ref{rem:ug}, we can calculate the Dirichlet-to-Neumann map
using \emph{any} solution to the boundary value problem, so we choose
$c_i = 0$.  It follows immediately that
\[
  \frac{\p u_i}{\p\nu_i}(\theta_{i-1}) = \frac{\p u_i}{\p\nu_i}(\theta_{i}) = 0
\]
for each $i$, hence the two-sided normal derivative $\gN u$ vanishes on $\Gamma$, and 
\begin{equation}
	\DN(\Gamma,\lambda_*,\chi) h = \Pi_\chi(\gN u) = 0,
\end{equation}
as expected.

\appendix

\section{Weights, cuts and pair compatibility conditions}
\label{sec:weights}

In this section we elaborate on some of our constructions and their connection to previous literature. In Section~\ref{sec:PCC} we discuss the relationship between the  strong pair compatibility condition in Definition~\ref{def:SPCC1}, and the  weak pair compatibility that appeared in earlier works, such as \cite{HH:2005a}, where it was simply referred to as the \emph{pair compatibility condition}. In Section~\ref{ssec:weight} we describe the cutting construction of \cite{HeSu}, and  explain how it is related to the valid weights introduced in Definition~\ref{def:weight}.

\begin{figure}[hbt]
  \centering
  \subcaptionbox{\label{cutone}}{\includegraphics[page=1]{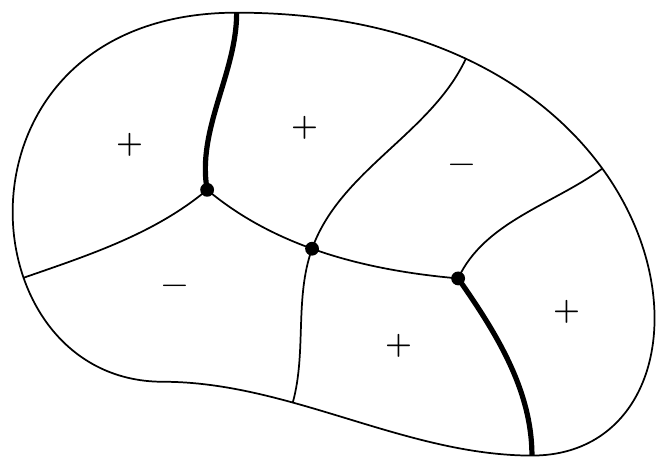}}
  \hskip 1cm
  \subcaptionbox{\label{cuttwo}}{\includegraphics[page=2]{figures2.pdf}}
  \\
  \subcaptionbox{\label{signsone}}{\includegraphics[page=4]{figures2.pdf}}
  \hskip 1cm
  \subcaptionbox{\label{signstwo}}{\includegraphics[page=5]{figures2.pdf}}
  \caption{Two valid cuts of the same partition are shown in (a) and (b), with the thick line denoting $\Gamma^*$ and the $\pm$ signs indicating the chosen orientations of each $D_i$. 
In (c) and (d) we show possible choices of $\{\chi_i\}$ for each of these cuts.
}
\label{fig:equivalentcuts}
\end{figure}

\subsection{Weak vs strong pair compatibility conditions}
\label{sec:PCC}

The strong pair compatibility condition (SPCC) was already described in Definition~\ref{def:SPCC1}, which we repeat here for convenience.

\begin{definition}
\label{def:SPCC2}
  A two-sided, weakly regular partition $\cD$ is said to satisfy the
  \emph{strong pair compatibility condition} (SPCC) if there exists a
  choice of positive ground states $\{u_i\}_{i=1}^k$
  for the Dirichlet Laplacians on $D_i$ such that, for any pair of neighbors
  $ D_i, D_j$, the function $u_{ij}$ defined by
  \begin{equation}
    u_{ij} \big|_{D_i} = u_i,
    \qquad
    u_{ij} \big|_{D_j} = -u_j,    
  \end{equation}
  is an eigenfunction of the Dirichlet Laplacian on
  $\Int (\overline {D_i\cup D_j})$.
\end{definition}

Nodal partitions obviously satisfy the SPCC. The same is true of spectral minimal partitions (see~\cite{HHOT}), and in Proposition~\ref{prop:chinodal_spcc} we showed that a partition 
satisfies the SPCC if and only if it is $\chi$-nodal. A partition satisfying 
the SPCC is necessarily an \emph{equipartition}, in the sense that the ground state energy (the smallest 
eigenvalue of the Dirichlet Laplacian) on each $D_i$ is the same. We denote this common value by $\lambda(\cD)$.

We next recall the weak pair compatibility condition.

\begin{definition}
\label{def:WPCC}
  A two-sided, weakly regular equipartition $\cD$ is said to satisfy the
  \emph{weak pair compatibility condition} (WPCC) if for each 
  pair of neighbors $ D_i, D_j$, there exists an eigenfunction of the 
  Dirichlet Laplacian on  $\Int (\overline {D_i\cup D_j})$ with eigenvalue $\lambda(\cD)$ 
  and nodal set $\p D_i \cap \p D_j$.
\end{definition}

\begin{remark}
By~\cite[Theorem~2.6]{HHOT} applied to each pair of neighbors, it follows that partitions that satisfy the WPCC also have the equal angle property; cf. Corollary~\ref{cor:equalangle}.
\end{remark}

It is obvious that SPCC implies WPCC. When $\Omega$ is simply connected, a bipartite equipartition satisfying WPCC is nodal, and hence satisfies SPCC, by \cite[Theorem~1.3]{HH:2005a}. If $\Omega$ is not simply connected, however, it is possible to find an equipartition (for a Schr\"odinger operator with $C^\infty$ potential) that satisfies WPCC but not SPCC, as shown in \cite[Section~7]{HH:2005a}.

\subsection{Weights and cuts}
\label{ssec:weight}
Assuming throughout that $\cD$ is a two-sided, weakly regular partition,
with nodal set $\Gamma$, we first decompose the smooth part of
$\Gamma$ into disjoint open curves, labeled $\{C_a\}$, so that
$\Gamma = \overline{\cup_a C_a}$. Since $\cD$ is two-sided, each $C_a$ is contained in
$\Gamma_i \cap \Gamma_j$ for some $i \neq j$. Without loss of
generality we can assume $i < j$, and we denote these labels by $i(a)$
and $j(a)$.

\begin{definition}
  \label{def:valid_cut}
  A subset $\cC \subset \{C_a\}$ is called a \emph{valid cut} of the
  partition $\cD$ if there exists a choice of orientations on the
  subdomains $\{D_i\}$ such that $C_a \in \cC$ if and only if
  $D_{i(a)}$ and $D_{j(a)}$ have the same orientation.
\end{definition}

It is sometimes convenient to identity a subset $\cC = \{C_{a_1}, \ldots, C_{a_p}\} \subset \{C_a\}$ with the corresponding closed subset
\begin{equation}
	\Gamma^* := \overline{C_{a_1} \cup \cdots \cup C_{a_p}}
\end{equation}
of $\Gamma$.
We mention that $\cC \subset \{C_a\}$ is a valid cut if $\Gamma 
\setminus \Gamma^*$ is a $\mathbb{Z}_2$-homological 1-cycle of $\Omega$ (viewed
as a cell complex) relative to the boundary $\partial \Omega$.  It is
immediate that the empty set is a valid cut of $\cD$ if and only if
$\cD$ is bipartite.

The \emph{maximal cut} $\cC = \{C_a\}$, for which $\Gamma^* = \Gamma$, is always valid\,---\,it corresponds
to all subdomains having the same orientation.  However, usually one
is interested in cuts that are as small as possible. We thus say that a cut is \emph{minimal} if $\Omega \setminus \Gamma^*$ is connected.

\begin{proposition}\cite[Prop~4.2]{HeSu}
  \label{prop:spcc-nodal2}
There exists a minimal valid cut $\cC \subset \{C_a\}$.
\end{proposition}

Finally, we describe how valid cuts are related to the valid weights $\{\chi_i\}$ in Definition~\ref{def:weight}. Given a set of valid weights $\{\chi_i\}$, we obtain a valid cut $\cC$ 
by declaring that $C_a \in \cC$ if and only if $\chi_{i(a)} = - \chi_{j(a)}$. That is, the cut set $\Gamma^*$ is the union of all $\Gamma_i \cap \Gamma_j$ along which $\chi_i = - \chi_j$. More 
precisely, we have the following.

\begin{proposition}
  \label{prop:valid_cuts_valid_weights}
Valid cuts are in one-to-one correspondence with edge-equivalence classes of valid weights.
\end{proposition}

\begin{proof}
Given a valid cut, i.e.\ a choice of orientation for each $D_i$, we get an induced orientation on each $\p D_i$. Choosing an orientation on each smooth component of $\Gamma$, we obtain a valid set of weights $\{\chi_i\}$ with the property that $\chi_i = -\chi_j$ if and only if $\p D_i \cap \p D_j$ is in the cut set $\Gamma^*$. Changing the orientation on any smooth part of $\Gamma$ will give a different, but edge equivalent, set of weights (recall Definition~\ref{def:weight_equivalences}), so we get a map from valid cuts to edge-equivalence 
classes of valid weights. Conversely, a set of valid weights gives an orientation on each $D_i$, and hence a valid cut. It is easily seen that edge-equivalent weights generate the same cut.
\end{proof}

\begin{remark}
  \label{rem:orientation_independent_valid}
  The proof of Proposition~\ref{prop:valid_cuts_valid_weights}
  suggests an equivalent way to define valid cuts and weights: a cut
  $\Gamma^*$ is valid if a generic closed path in $\Omega$
  intersects $\Gamma \setminus \Gamma^*$ an even number of times, and a
  choice of weights $\{\chi_i\}$ is valid if the set
  $\{C_a : \chi_{i(a)} = - \chi_{j(a)}\}$
defines a valid cut.
  This alternative definition is not as constructive as
  Definition~\ref{def:weight}, but it has the advantage of not depending on the 
	manifold structure of $\Omega$, and
  is thus more convenient for considering
  partitions on metric graphs.
\end{remark}

\begin{remark}
Another way of viewing the constructions in this paper is to introduce Aharonov--Bohm operators, as in~\cite{HeSu}. Given a set of weights $\chi$ that generates a minimal valid cut, the corresponding $\Delta^\chi$ is equivalent to a certain Aharonov--Bohm operator, with Aharonov--Bohm solenoids with flux $\pi$ placed at the singular
  points $x_\ell$ of $\Gamma$ for which $\nu_\ell$ is odd (recall Definition~\ref{def:regular}).
\end{remark}

\section{Explicit construction of the canonical solution to \eqref{eq:Di2}}
\label{sec:equivalent}

In this section we give an alternate, more explicit proof of the second claim in Theorem~\ref{thm:Aform}, regarding the existence of a \enquote{canonical solution} $\tilde u \in \dom(t^\chi)$ such that $\gN \tilde u \in S_\chi$ and $\tilde u_i$ solves~\eqref{eq:Di2} for each $i$. To 
do this we write the condition $\gN \tilde u \in S_\chi$ as a finite system of linear equations and then, by analyzing the corresponding matrix, prove that a solution always exists.

Fix $g \in \dom(A)$. For each $i$, the general solution of~\eqref{eq:Di2} is given by
\begin{equation}\label{eq:ui}
  u_i=u_i^g+c_i\phi_{*,i}
\end{equation}
for some $c_i \in \bbR$. Since $g \in \dom(A)$, we know from \eqref{Adom} that the two-sided normal derivative $\gN u$ is a function in $L^2(\Gamma)$, and is given by 
$\chi_i \partial_{\nu_i}u_i + \chi_j\partial_{\nu_j}u_j$ on $\Gamma_i\cap\Gamma_j$. This will be an element of the subspace $S_\chi$ if and only if
\begin{equation}
\label{cond1}
	I_i := \int_{\Gamma_i} \chi_i (\gN u)\partial_{\nu_i}\phi_{*,i} = 0
\end{equation}
for each $i$. Since each point in the smooth part of $\Gamma_i$ is contained in precisely one other $\Gamma_j$, we can rewrite this integral as
\begin{align}
\begin{split}\label{eq:system}
	I_i
	&= \sum_{j \neq i} \int_{\Gamma_i\cap\Gamma_j} \big(\partial_{\nu_i}u_i +\chi_{ij} \partial_{\nu_j}u_j\big)\,\partial_{\nu_i}\phi_{*,i} \\
	&= \sum_{j \neq i} \int_{\Gamma_i\cap\Gamma_j} \big(\partial_{\nu_i}u^g_i +\chi_{ij} \partial_{\nu_j}u^g_j\big)\,\partial_{\nu_i}\phi_{*,i}
	 +  \sum_{j \neq i} \int_{\Gamma_i\cap\Gamma_j} \big(c_i\partial_{\nu_i}\phi_{*,i} + \chi_{ij} c_j \partial_{\nu_j}\phi_{*,j}\big)\,\partial_{\nu_i}\phi_{*,i},
\end{split}
\end{align}
where we have denoted $\chi_{ij} = \chi_i \chi_j$ for convenience.
Let us introduce the notations
\[
	\alpha_{i,i} = 0, \qquad \alpha_{i,j} = \int_{\Gamma_i\cap\Gamma_j}|\partial_{\nu_i}\phi_{*,i}|^2, \quad i \neq j.
\]
It follows from~\eqref{eq:transmission} that $|\partial_{\nu_i}\phi_{*,i}| = |\partial_{\nu_j}\phi_{*,j}|$ on $\Gamma_i \cap \Gamma_j$, and so $\alpha_{i,j} = \alpha_{j,i}$ for all $i,j$. We similarly get
\[
	\int_{\Gamma_i \cap \Gamma_j} \chi_{ij} (\p_{\nu_j} \phi_{*,j})(\p_{\nu_i} \phi_{*,i}) 
	= -\alpha_{i,j},
\]
We then define
\begin{equation}\label{eq:di}
	d_i
	= - \sum_{j \neq i} \int_{\Gamma_i\cap\Gamma_j} \big(\partial_{\nu_i}u^g_i +\chi_{ij} \partial_{\nu_j}u^g_j\big)\,\partial_{\nu_i}\phi_{*,i}
\end{equation}
so the equation~\eqref{eq:system} becomes
\begin{equation}
\label{eq:equations}
	\sum_{j \neq i}(c_i - c_j) \alpha_{i,j} = d_i.
\end{equation}
We write the resulting system of equations in matrix form as
\begin{equation}\label{eq:thecondition}
  \underbrace{\begin{bmatrix}
           \sum_j \alpha_{1,j} & -\alpha_{1,2} & \dots & -\alpha_{1,k}\\
        -\alpha_{2,1} & \sum_j \alpha_{2,j} & \dots & -\alpha_{2,k}\\
        \vdots & \vdots &\ddots & \vdots\\
        -\alpha_{k,1} & -\alpha_{k,2} & \dots & \sum_j \alpha_{k,j}\\
  \end{bmatrix}}_{=:A}
  \begin{bmatrix}
        c_1\\c_2\\ \vdots \\ c_k\\
  \end{bmatrix}
  =
  \begin{bmatrix}
    d_1\\d_2\\ \vdots \\ d_k\\
  \end{bmatrix},
\end{equation}
and observe that the vector $(c_1,c_2,\ldots,c_k)^t=(1,1,\ldots,1)^t$ lies in the kernel of the matrix $A$. 

Without loss of generality, we can label the domains $\{D_i\}$ in the partition inductively so that $D_{i+1}$ is a neighbor of at least one of $D_1, \ldots, D_i$, with $D_1$ arbitrary. For the numbers $\alpha_{i,j}$, this means that
\begin{equation}\label{eq:neighbor}
\begin{gathered}
\alpha_{1,2}>0,\\
\alpha_{1,3}+\alpha_{2,3}>0,\\
\alpha_{1,4}+\alpha_{2,4}+\alpha_{3,4}>0,\\
\vdots\\
\alpha_{1,k}+\alpha_{2,k}+\cdots+\alpha_{k-1,k}>0.
\end{gathered}
\end{equation}

\begin{lemma}
Let $A$ be the symmetric $k\times k$ matrix in~\eqref{eq:thecondition} and assume that the inequalities in~\eqref{eq:neighbor} hold. Then $\ker A$ is spanned by $(1,1,\ldots,1)^t$.
\end{lemma}
    
\begin{proof}
Consider the quadratic form $q[\mathbf c]=\langle A\mathbf c,\mathbf c\rangle$ corresponding to the matrix $A$ above, where $\mathbf c=(c_1,c_2,\ldots,c_k)^t$. From~\eqref{eq:equations} we find that the quadratic form $q[\mathbf c]$ can be written as
\[
\begin{gathered}
  \alpha_{1,2}(c_1-c_2)^2\\
  +\alpha_{1,3}(c_1-c_3)^2+\alpha_{2,3}(c_2-c_3)^2\\
  +\alpha_{1,4}(c_1-c_4)^2+\alpha_{2,4}(c_2-c_4)^2+\alpha_{3,4}(c_3-c_4)^2\\
  +\cdots\\
  +\alpha_{1,k}(c_1-c_k)^2+\alpha_{2,k}(c_2-c_k)^2+\cdots+\alpha_{k-1,k}(c_{k-1}-c_k)^2.
\end{gathered}   
\]
Since $\alpha_{i,j} \geq 0$ for all $i,j$, we see that $q$ (and hence $A$) is non-negative. It remains to identify the kernel. Assume that $q[\mathbf c]=0$ for some $\mathbf c$. Then, reading from the top line above, we conclude that $c_2=c_1$ since $\alpha_{1,2}>0$. Inserting $c_1=c_2$, we conclude from the next row that $c_3=c_2$ since $\alpha_{1,3}+\alpha_{2,3}>0$. Continuing in this manner, we conclude that $c_k=c_{k-1}=\cdots=c_2=c_1$. This means that the kernel of $A$ is spanned by the vector $(1,1,\ldots,1)^t$.
\end{proof}

Finally, from~\eqref{eq:di} we observe that $\sum d_i$ contains a term
\[
	- \int_{\Gamma_i\cap\Gamma_j} \left\{ \big(\partial_{\nu_i}u^g_i +\chi_{ij} \partial_{\nu_j}u^g_j\big)\,\partial_{\nu_i}\phi_{*,i} + \big(\partial_{\nu_j}u^g_j +\chi_{ij} \partial_{\nu_i}u^g_i\big)\,\partial_{\nu_j}\phi_{*,j} \right\}
\]
for each pair of neighboring domains, and by \eqref{eq:transmission} each of the integrands vanishes. This means $\sum d_i = 0$, so the vector $\mathbf{d}=(d_1,d_2,\ldots,d_k)^t$ is orthogonal to the kernel of $A$. Thus, the system~\eqref{eq:thecondition} will always be solvable, by the Fredholm alternative for symmetric matrices.

\bibliographystyle{plain}
\bibliography{nodal}

\end{document}